%
%
%
%
%
\hsize=5in
\baselineskip=12pt
\vsize=20.4cm
\parindent=.5cm
\predisplaypenalty=0
\hfuzz=2pt
\frenchspacing
\def\latexfmt{latex}
\ifx\fmtname\latexfmt\else
%
%
\input amssym.def
\font\ss=cmss10
\def\titlefonts{\baselineskip=1.44\baselineskip
	\font\titlef=cmbx12
	\font\titlei=cmmib10 scaled 1200
	\skewchar\titlei='177
	\textfont1=\titlei
	\titlef
	}
\font\tenib=cmmib10 
\font\sevenib=cmmib7
\font\fiveib=cmmib5
\skewchar\tenib='177
\skewchar\sevenib='177
\skewchar\fiveib='177
\font\tenbsy=cmbsy10 
\font\sevenbsy=cmbsy7
\font\fivebsy=cmbsy5
\skewchar\tenbsy='60
\skewchar\sevenbsy='60
\skewchar\fivebsy='60
\def\boldfonts{\bf
	\textfont0=\tenbf
	\scriptfont0=\sevenbf
	\scriptscriptfont0=\fivebf
	\textfont1=\tenib
	\scriptfont1=\sevenib
	\scriptscriptfont1=\fiveib
	\textfont2=\tenbsy
	\scriptfont2=\sevenbsy
	\scriptscriptfont2=\fivebsy
	}
\font\ninerm=cmr9
\font\ninebf=cmbx9
\font\ninei=cmmi9
\skewchar\ninei='177
\font\ninesy=cmsy9
\skewchar\ninesy='60
\font\nineit=cmti9
\def\reffonts{\baselineskip=0.9\baselineskip
	\textfont0=\ninerm
	\def\rm{\fam0\ninerm}%
	\textfont1=\ninei
	\textfont2=\ninesy
	\textfont\bffam=\ninebf
	\def\bf{\fam\bffam\ninebf}%
	\def\it{\nineit}%
	}
%
%
\def\frontmatter{\vbox{}\vskip1cm\bgroup
	\leftskip=0pt plus1fil\rightskip=0pt plus1fil
	\parindent=0pt
	\parfillskip=0pt
	\pretolerance=10000
	}
\def\endfrontmatter{\egroup\bigskip}
\def\title#1{{\titlefonts#1\par}}
\def\author#1{\bigskip#1\par}
\def\address#1{\bigskip{\reffonts\it#1}}
\def\email#1{\bigskip{\reffonts{\it E-mail: }\rm#1}}
\def\thanks#1{\footnote{}{\reffonts\rm\noindent#1\hfil}}
\def\section#1\par{\ifdim\lastskip<\bigskipamount\removelastskip\fi
	\penalty-250\bigskip
	\vbox{\leftskip=0pt plus1fil\rightskip=0pt plus1fil
	\parindent=0pt
	\parfillskip=0pt
	\pretolerance=10000{\boldfonts#1}}\nobreak\medskip
	}

\fi
\def\proclaim#1. {\medbreak\bgroup{\noindent\bf#1.}\ \it}
\def\endproclaim{\egroup
	\ifdim\lastskip<\medskipamount\removelastskip\medskip\fi}
\def\item#1 #2\par{\ifdim\lastskip<\smallskipamount\removelastskip\smallskip\fi
	{\rm#1}\ #2\par\smallskip}
\def\Proof#1. {\ifdim\lastskip<\medskipamount\removelastskip\medskip\fi
	{\noindent\it Proof.}\ }
\def\endproof{\quad\hbox{$\square$}\medskip}
\def\Remark. {\ifdim\lastskip<\medskipamount\removelastskip\medskip\fi
	{\noindent\bf Remark.}\quad}
\def\Remarks. {\ifdim\lastskip<\medskipamount\removelastskip\medskip\fi
	{\noindent\bf Remarks.}\quad}
\def\endremark{\medskip}
\def\Example#1. {\ifdim\lastskip<\medskipamount\removelastskip\medskip\fi
	{\bf Example #1.} }
%
%
\newcount\citation
\newtoks\citetoks
\def\citedef#1\endcitedef{\citetoks={#1\endcitedef}}
\def\endcitedef#1\endcitedef{}
\def\citenum#1{\citation=0\def\curcite{#1}%
	\expandafter\checkendcite\the\citetoks}
\def\checkendcite#1{\ifx\endcitedef#1?\else
	\expandafter\lookcite\expandafter#1\fi}
\def\lookcite#1 {\advance\citation by1\def\auxcite{#1}%
	\ifx\auxcite\curcite\the\citation\expandafter\endcitedef\else
	\expandafter\checkendcite\fi}
\def\cite#1{\makecite#1,\cite}
\def\makecite#1,#2{[\citenum{#1}\ifx\cite#2]\else\expandafter\clearcite\expandafter#2\fi}
\def\clearcite#1,\cite{, #1]}
%
%
\def\references{\section References\par
	\bgroup
	\parindent=0pt
	\reffonts
	\rm
	\frenchspacing
	\setbox0\hbox{99. }\leftskip=\wd0
	}
\def\endreferences{\egroup}
\newtoks\nextauth
\newif\iffirstauth
\def\checkendauth#1{\ifx\endauth#1%
		\iffirstauth\the\nextauth
		\else{} and \the\nextauth\fi,
	\else\iffirstauth\the\nextauth\firstauthfalse
		\else, \the\nextauth\fi
		\expandafter\auth\expandafter#1\fi
	}
\def\auth#1,#2;{\nextauth={#1 #2}\checkendauth}
\newif\ifinbook
\newif\ifbookref
\def\nextref#1 {\par\hskip-\leftskip
	\hbox to\leftskip{\hfil\citenum{#1}.\ }%
	\initnextref}
\def\initnextref{\bookreffalse\inbookfalse\firstauthtrue\ignorespaces}
\def\paper#1{{\it#1,}}

\def\book#1{\bookreftrue{\it#1},}
\def\journal#1{#1\ifinbook,\fi}
\def\bookseries#1{#1,}
\def\Vol#1{\ifbookref Vol. #1,\else\ifinbook Vol. #1,\else{\bf#1}\fi\fi}
\def\nombre#1{no.~#1}
\def\publisher#1{#1,}
\def\Year#1{\ifbookref #1.\else\ifinbook #1,\else(#1)\fi\fi}
\def\Pages#1{\ifinbook pp. #1.\else #1.\fi}

%
%
\newsymbol\square 1003
\mathcode`\#="2023
\let\<\langle
\let\>\rangle
\def\*#1{\vphantom{#1}^*#1}
\def\ad{\mathop{\fam0 ad}\nolimits}
\def\alg{\hbox{\rm-\ss alg}}
\def\chr{\mathop{\fam0 char}\nolimits}
\def\coalg{\hbox{\rm-\ss coalg}}
\def\cop{^{\mathop{\fam0 cop}}}
\def\End{\mathop{\fam0 End}\nolimits}
\def\Hom{\mathop{\fam0 Hom}\nolimits}
\def\id{\mathop{\fam0 id}\nolimits}
\def\Im{\mathop{\fam0 Im}}
\def\Lie{\mathop{\fam0 Lie}}
\def\Mat{\mathop{\fam0 Mat}\nolimits}
\def\mo{\hbox{\rm-\ss mod}}
\def\mod#1{\ifinner\mskip8mu(\mathop{\rm mod}#1)
        \else\mskip12mu(\mathop{\rm mod}#1)\fi}
\def\op{^{\mathop{\fam0 op}}}
\def\rk{\mathop{\fam0 rk}\nolimits}
\def\tr{\mathop{\fam0 tr}\nolimits}
\def\textsum{\textstyle\sum\limits}
\def\iso{\mathrel{\setbox0\hbox{$\rightarrow$}%
  \setbox1\hbox{$\sim$}%
  \dimen0=\wd0\advance\dimen0by-\wd1\advance\dimen0by-.05em
  \box0\kern-\wd1\kern-\dimen0\raise0.75ex\box1\kern\dimen0}}
\def\mapright#1{{}\mathrel{\mathop{\longrightarrow}\limits^{#1}}{}}
\def\mapleft#1{{}\mathrel{\mathop{\longleftarrow}\limits^{#1}}{}}
\def\longmapright#1#2{{}\mathrel{\smash{\mathop{\count0=#1
  \loop
    \ifnum\count0>0
    \advance\count0 by-1\mathord-\mkern-4mu
  \repeat
  \mathord\rightarrow}\limits^{#2}}}{}}
\def\mapdown#1#2{\llap{$\vcenter{\hbox{$\scriptstyle{#1}$}}$}\big\downarrow
  \rlap{$\vcenter{\hbox{$\scriptstyle{#2}$}}$}}
\def\mapup#1#2{\llap{$\vcenter{\hbox{$\scriptstyle{#1}$}}$}\big\uparrow
  \rlap{$\vcenter{\hbox{$\scriptstyle{#2}$}}$}}
\def\diagram#1{\vbox{\halign{&\hfil$##$\hfil\cr #1}}}
\def\diagramskip{\noalign{\smallskip}}
\def\mlines$$#1$$#2\mlines{{\belowdisplayskip=0pt
  \belowdisplayshortskip=\belowdisplayskip$$#1$$
  \abovedisplayshortskip=\the\medskipamount
  #2}\vskip\belowdisplayskip\noindent}
\newcount\tablecolumns
\newcount\tcolcount
\def\checkendtable#1{\ifx\endtable#1\noalign{\hrule}\egroup\egroup
	\else\tablerow#1\fi}
\def\endtable{endtable}
\def\table#1 {\tablecolumns#1\vbox\bgroup\offinterlineskip
	\halign\bgroup\vrule##&&\hfil\hskip4pt$##$\hskip4pt\hfil\vrule\cr
	\checkendtable}
\def\tablerow#1\cr{\noalign{\hrule}%
	height2pt\global\tcolcount=\tablecolumns&\emptyentry
	\strut&#1\cr
	height2pt\global\tcolcount=\tablecolumns&\emptyentry
	\checkendtable}
\def\emptyentry{\ifnum\the\tcolcount=1\cr
	\else\global\advance\tcolcount by-1&\emptyentry\fi}
\let\al\alpha
\let\be\beta
\let\ga\gamma
\let\de\delta
\let\ep\varepsilon
\let\io\iota
\let\ka\kappa
\let\la\lambda
\let\om\omega
\let\ph\varphi
\let\si\sigma
\let\th\theta
\let\De\Delta
\let\Up\Upsilon
\def\c{\tilde c}
\def\g{{\frak g}}
\def\pit{\widetilde\pi}
\def\r{\overline r}
\def\rhot{\overline\rho}
\def\R{\widetilde R}
\def\S{{\frak S}}
\citedef
Ali01
And01
Bea
Bl89
Chase69
Coh90
Doi86
Doi89
Dr87
Dr89a
Dr89b
Et00
Fai
Gr96
Jac
Jan
Joy93
Kas
Kop
Kr76
Kr81
Lar69
Maj94
Maj
Man
Mas
Mil75
Mo
Nich89
Oo74
Oo76
Oo80
Pre99
Rad76
Rad77
Rad92
Rad93
Scha96
Scha97
Schn90
St78
St
Sw
Taft75
Taft80
Ul82
Oy94
Wei71
\endcitedef
%
%
\frontmatter

\title{Hopf Galois extensions, triangular structures,\break
and Frobenius Lie algebras in prime characteristic}
\author{Serge Skryabin}
\address{Chebotarev Research Institute,
Kazan, Russia}
\address{Current address: Mathematisches Seminar, University of Hamburg,\break
Bundesstr. 55, 20146 Hamburg, Germany}
\email{fm1a009@math.uni-hamburg.de}
\thanks{The author gratefully acknowledges the support of the Alexander
von Humboldt Foundation through the program of long term cooperation}

\endfrontmatter

\section
Introduction

The final goal of this paper is to introduce certain finite dimensional Hopf
algebras associated with restricted Frobenius Lie algebras over a field of characteristic
$p>0$. The antipodes of these Hopf algebras have order either $2p$ or 2, and
in the minimal dimension $p^2$ there exists just one Hopf algebra in this
class which coincides with an example due to Radford \cite{Rad77} of a Hopf
algebra with a nonsemisimple antipode (the first example of this kind in
\cite{Taft75} does not quite fit into our scheme). Another feature of the Hopf
algebras under consideration is that they admit triangular structures of
maximal rank, so that they are isomorphic to their dual Hopf algebras taken
with the opposite multiplications. Regarding the algebra structure alone,
these Hopf algebras are isomorphic to the restricted enveloping algebras of
the corresponding Lie algebras (such isomorphisms are not canonical). In particular,
they are never semisimple. It should be mentioned that all semisimple cosemisimple
triangular Hopf algebras over any algebraically closed field were classified by Etingof
and Gelaki \cite{Et00}.

The Hopf algebras we deal with are instances of transformations in Hopf Galois extensions
discovered not so long ago. The Hopf Galois theory originated in the work of Chase and
Sweedler \cite{Chase69} and received a full treatment by Kreimer and Takeuchi \cite{Kr81}.
We recall the notion of a Hopf Galois algebra in section 1 of the paper. Suppose that $H$ is a
finite dimensional Hopf algebra over a field and $A$ is a left $H$-module
algebra which is $H^*$-Galois with respect to the corresponding right
$H^*$-comodule structure. A construction of Schauenburg \cite{Scha96}
generalizing a special case of commutative $A$ and $H^*$ considered by Van
Oystaeyen and Zhang \cite{Oy94} produces a certain Hopf algebra $L(A,H^*)$. As
was pointed out by Greither \cite{Gr96} in the commutative case the dual Hopf
algebra is isomorphic to the endomorphism algebra $E=\End_HA$. A similar description is
valid in general. In section 1 we show that the Hopf algebra structure on
$E$ can be easily introduced directly using properties of a duality functor
$D$ between the categories of finite dimensional $H$-modules and $E$-modules.
It is this Hopf algebra $E$ that we are interested in. Another way to construct
$E$ is to twist the comultiplication in $H$ by means of Drinfeld's process
\cite{Dr89b}. Such a twist is determined by an element of $H\otimes H$ which is
unique, however, only up to a ``gauge equivalence", and in the situation of
primary interest to us it is hardly possible to see such an element. The Hopf
Galois theory provides a way to handle the subject properly.

There is a bijective correspondence between the quasitriangular structures on $H$ and $E$.
In section 2 we investigate a relationship between the ranks of two mutually corresponding
quasitriangular structures $R$ on $H$ and $\R$ on $E$. It is shown in Theorem 2.5 that $\dim A$ divides
the product $(\rk R)(\rk\R)$ provided that $A$ is central simple. In particular, $\rk\R=\dim A=\dim E$,
that is, $\R$ has maximal rank when $H$ is
cocommutative and $R=1\otimes1$. In section 3 we show that the squares of antipodes in both $H$
and $E$ are closely related to a certain automorphism $\th$ of the algebra $A$. Theorem 3.8 determines
the order of the antipode of $E$ in case when $H$ is cocommutative and $A$ is central simple.
It is described in terms of the order of the modular function $\al\in H^*$ in the group of invertible
elements of $H^*$. Section 4 explains how each $H$-invariant, with respect to
a suitable module structure, generating subspace in the algebra $H$ or $A$ gives rise to a generating
subspace in $E$ or $E^*$, respectively. If $H$ is cocommutative and $A$ is
central simple, then $E^*$ is mapped onto $E$, giving a second kind of
generating subspaces in $E$. Particular elements occurring in such generating
subspaces are introduced and their properties are recorded.

All the results discussed above are obtained in general Hopf theoretical
settings. Since the algebras over a field are of main concern in this paper,
basic facts from Hopf Galois theory are recalled in much less generality than
they are available in the literature. Another preference of using the module
structures rather than the comodule ones stems from the desire to find a link
between the Hopf algebras and the representation theory, especially in case of Lie
algebras.

Suppose now that $\g$ is a finite dimensional $p$-Lie algebra over a field $k$
such that $\chr k=p$. To each $\xi\in\g^*$ there corresponds a finite
dimensional factor algebra $U_\xi(\g)$ of the universal enveloping algebra of $\g$
\cite{St, 5.3}. Here $U_0(\g)$ is a cocommutative Hopf algebra in a natural
way, and each $U_\xi(\g)$ is a left $U_0(\g)$-module algebra with respect to
the adjoint action. It turns out that $U_\xi(\g)$ is $U_0(\g)^*$-Galois if
and only if $U_\xi(\g)$ is central simple. In this case one obtains a Hopf
algebra $E_\xi=\End_\g U_\xi(\g)$ as was explained above. In section 6
preceding results are reformulated for $E_\xi$. There is a canonical way to
pick out a generating subspace for $E_\xi$ and a multiplication on this
subspace which defines a structure not of an ordinary Lie algebra but rather
of a quantum modification. What one gets here is covered by the concept of
Lie algebras in symmetric tensor categories proposed by Manin \cite{Man}. In
this sense $E_\xi$ is an enveloping algebra of a quantum Lie algebra.

The $p$-Lie algebra $\g$ is called {\it Frobenius} if there exists $\xi\in\g^*$ such that
the associated alternating bilinear form $\be_\xi(x,y)=\xi([xy])$ defined on $\g$ is
nondegenerate \cite{Oo74}. If $\g$ is Frobenius and $k$ is algebraically closed then it
follows from the results in \cite{Mil75} and \cite{Pre99} that the family of reduced
enveloping algebras $U_\xi(\g)$ contains a simple algebra. It is an open problem if the
converse is true, and such a result would parallel the characterization due to
Ooms \cite{Oo74}, \cite{Oo76} of finite dimensional Lie algebras over a field of characteristic 0
having primitive universal enveloping algebras. Explicit computations of Hopf
algebras $E_\xi$ for low dimensional Frobenius $\g$ are done at the end of the paper.

\section
1. Transformation of a Hopf algebra

We work over the ground field $k$. The indication of $k$ in $\otimes$, $\Hom$,
$\End$ will be omitted. Let $H$ be a finite dimensional Hopf algebra over $k$.
All unexplained notions related to Hopf algebras are taken from \cite{Mo} and
\cite{Sw}. Denote by $\De$, $\ep$, $S$ the comultiplication, counity and
antipode in $H$ and other Hopf algebras. We write symbolically $\De
h=\sum_{(h)}h'\otimes h''$ for $h\in H$. Denote by $H\mo$ the category of finite
dimensional $H$-modules. Given $U,V\in H\mo$, the action of $H$ in
$U\otimes V$ will be always assumed to come via $\De$. There are two
ways to make $H$ operate in the dual of $V$: one takes $\<h\xi,v\>$ to equal
either $\<\xi,S(h)v\>$ or $\<\xi,S^{-1}(h)v\>$ for $\xi\in V^*$ and $v\in V$.
One thus obtains the left and right duals of $V$ in $H\mo$, denoted as $V^*$ and
$\*V$, respectively. Put $V^H=\{v\in V\mid hv=\ep(h)v$ for all $h\in H\}$.
For $U,V,W$ in $H\mo$ there are linear bijections
\mlines
$$
(\*V\otimes W)^H\cong\Hom_H(V,W)\cong(W\otimes V^*)^H,\eqno(1.1)
$$$$
\Hom_H(V,\*U\otimes W)\cong\Hom_H(U\otimes V,W)\cong\Hom_H(U,W\otimes V^*)\eqno(1.2)
$$
\mlines
(see \cite{Kas, Prop.~XIV.2.2}). In particular, if $\mu:U\otimes V\to W$ is
a morphism in $H\mo$, then the set $\{v\in V\mid\mu(U\otimes v)=0\}$ is an
$H$-submodule of $V$. Indeed, it coincides with kernel of the $H$-module map
$V\to\*U\otimes W$ corresponding to $\mu$.
When $k$ stands for an object of $H\mo$, the action of $H$ is given by $\ep$.
By a similar convention $H$ operates in $H\in H\mo$ via left multiplications.
On the other hand, $H_{\ad}$ presupposes the adjoint action under which
$h\triangleright g=\sum_{(h)}h'gS(h'')$ for $h,g\in H$.
The Hopf algebra $H\op$ has the same comultiplication as in $H$ and the opposite
multiplication. In $H\cop$ the comultiplication is the opposite one:
$\De\op h=\sum_{(h)}h''\otimes h'$.

Let now $A$ be a left $H$-module algebra, so that the unity map $k\to A$ and
the multiplication map $A\otimes A\to A$ are $H$-module homomorphisms.
Changing the multiplication in $A$ to the opposite one produces a left
$H\cop$-module algebra $A\op$. Define linear maps $\pi,\pi':A\otimes H\to\End A$
and $\ga,\ga':A\otimes A\to A\otimes H^*\cong\Hom(H,A)$ by
\mlines
$$
\pi(a\otimes h)(b)=a(hb),\eqno(1.3)
$$$$
\pi'(a\otimes h)(b)=(hb)a,\eqno(1.4)
$$$$
\ga(a\otimes b)(h)=a(hb),\eqno(1.5)
$$$$
\ga'(a\otimes b)(h)=(ha)b,\eqno(1.6)
$$
\mlines
where $a,b\in A$ and $h\in H$. Two of these maps give homomorphisms of
algebras $\pi:A\#H\to\End A$ and $\pi':A\op\#H\cop\to\End A$ where
$\#$ denotes the smash product algebra structure.

\proclaim
Proposition 1.1.
If one of the maps $\pi,\pi',\ga,\ga'$ is bijective then so are all others.
\endproclaim

This is contained in \cite{Kr81} and \cite{Ul82}. The only difference is that
we dropped the assumption $A^H=k$ on the subalgebra of $H$-invariants $A^H\subset A$. The
map $\ga'$ is $A$-linear with respect to right multiplications on the second
tensorand in $A\otimes A$ and the first tensorand in $A\otimes H^*$.
One obtains $\pi$ from $\ga'$ by applying the functor $\Hom_A(?,A)$ after
natural identifications. Since the $A$-modules $A\otimes A$ and $A\otimes H^*$
are free, $\pi$ is bijective if and only if so is $\ga'$. There is a similar
relationship between $\pi'$ and $\ga$. Finally, $\ga'=\Phi\circ\ga$ for a
suitable invertible transformation $\Phi$ of $A\otimes H^*$ constructed in
\cite{Kr81, Prop.~1.2}.

One can regard $A$ as a right $H^*\!$-comodule algebra with respect to the
structure map $A\to A\otimes H^*$ such that $a\mapsto\ga(1\otimes a)$. The algebra $A$ is
called {\it right $H^*\!$-Galois\/} or an {\it$H^*\!$-Galois extension} of $k$ if
$\pi,\pi',\ga,\ga'$ are bijective \cite{Chase69}, \cite{Kr76}. It suffices
to assume only that $\ga$ is surjective provided that $A^H=k$ \cite{Kr81,
(1.4) and (1.7)}. If $A$ is $H^*\!$-Galois then $\dim A<\infty$. Left comodule Galois algebras are defined similarly.

\proclaim
Proposition 1.2.
Suppose that $A$ is $H^*\!$-Galois. Then:

\item(i)
As an $H$-module, $A$ is free of rank $1$. In particular, $A^H=k$.

\item(ii)
The algebra $A$ is Frobenius. If $\chi\in A^*$ spans a onedimensional
$H$-submodule, then the bilinear form $(a,b)=\chi(ab)$
defined on $A$ is nondegenerate.

\item(iii)
For any $A\#H$-module $M$ the linear map $A\otimes M^H\to M$ given by
$a\otimes v\mapsto av$ is bijective.

\endproclaim

\Proof.
Assertion (i) is the so-called normal basis property. It follows quickly
from the Krull-Schmidt Theorem and the fact that $A\#H$ is a free left
$H$-module of rank $d=\dim A$ which is isomorphic, on the other hand, to a
direct sum of $d$ copies of $A$ since $A\#H\cong\End A$. An
equivalent formulation is given in \cite{Kr76, Prop.~2}. Since $H$ is
a Frobenius algebra by \cite{Lar69}, one has $\dim A^H=1$.

Assertion (ii) is a special case of \cite{Kr81, Th. 1.7(5)}. By the assumption
on $\chi$ there exists an algebra homomorphism $\al:H\to k$ such that
$\chi(ha)=\al(h)\chi(a)$ for all $h\in H$ and $a\in A$. Then $\ker\chi$ is
stable under $H$, and it follows that so too is the left ideal $I=\{b\in
A\mid\chi(Ab)=0\}$ of $A$. As $I$ is a proper $A\#H$-submodule of $A$, it has to be
zero.

Assertion (iii) is a consequence of the fact that $M$ is a direct sum of copies of
the simple $A\#H$-module $A$. Category equivalences in more general settings are
described in \cite{Coh90}, \cite{Doi86}, \cite{Doi89}, \cite{Schn90},
\cite{Ul82}.
\endproof

\proclaim
Corollary 1.3.
Let $V\in H\mo$, and let $A$ operate in $M=\Hom(V,A)$ by the rule
$(a\si)(v)=a\si(v)$ where $a\in A$, $\si\in M$ and $v\in V$.
The map $A\otimes\Hom_H(V,A)\to M$ afforded by this action of $A$ is bijective.
\endproclaim

\Proof.
Let $A$ operate in $A\otimes V^*$ via left multiplications on the first tensorand,
and let $H$ operate via $\De:H\to H\otimes H$ as usual. These two module structures
are compatible so that $A\otimes V^*$ becomes an $A\#H$-module. We make $M$ into an $A\#H$-module
using the canonical $k$-linear bijection $M\cong A\otimes V^*$.
According to (1.1) $M^H=\Hom_H(V,A)$. It remains to apply Prop.~1.2(iii).
\endproof

We now make a standing assumption for the most of the rest of the paper:
$$
\hbox{\it$A$ is assumed to be $H^*\!$-Galois, and $E=\End_HA$.}\eqno(1.7)
$$
By our convention both $H$ and $E$ operate on $A$ from the left. Since $A\cong
H$ in $H\mo$, there exist a noncanonical isomorphism of algebras
$E\cong H\op\cong H$ (the second isomorphism is given by the antipode $S$)
and an isomorphism of $E$-modules $A\cong E$.

\proclaim
Proposition 1.4.
Let $U,V\in H\mo$.

\item(i)
The contravariant functor $D=\Hom_H(?,A)$ gives a Morita duality
between $H\mo$ and $E\mo$. The inverse functor is $D^{-1}=\Hom_E(?,A)$.

\item(ii)
One has $\dim D(V)=\dim V$ for all $V$. In fact $D(V)\cong V^*$ under
noncanonical $k$-linear bijections which are natural in $V$.

\item(iii)
There are $k$-linear bijections $D(U)\otimes D(V)\iso D(U\otimes V)$,
$\,\tau\otimes\si\mapsto\tau\smile\si$, where
$(\tau\smile\si)(u\otimes v)=\tau(u)\si(v)$ for $\tau\in D(U)$,
$\si\in D(V)$, $u\in U$, $v\in V$.

\item(iv)
There are perfect pairings $D(V)\times D(V^*)\to k$ defined by the rule
$\<\si,\tau\>=\sum_i\si(e_i)\tau(e_i^*)$ for $\si\in D(V)$,
$\tau\in D(V^*)$ where $\{e_i\}$ and $\{e_i^*\}$
are dual bases of the vector spaces $V$ and $V^*$, respectively.

\item(v)
If $\{\si_i\}$ and $\{\tau_i\}$ are dual bases of $D(V)$ and $D(V^*)$
with respect to the pairing in {\rm(iv)}, then
$\<\xi,v\>=\sum_i\tau_i(\xi)\si_i(v)$ for all $\xi\in V^*$ and $v\in V$.

\endproclaim

\Proof.
Since $H$ is Frobenius, it follows from Proposition 1.2(i) that $A$ is an
injective cogenerator in $A\mo$. Hence $D$ is indeed a Morita duality
with $D^{-1}$ as stated \cite{Fai, Th.~23.25, $(7)\Rightarrow(8)$}.
Choosing a particular isomorphism $A\cong H\cong\*\!H$ in $H\mo$ and
using (1.2), we get bijections
$$
D(V)=\Hom_H(V,A)\cong\Hom_H(H\otimes V,k)\cong\Hom_H(H,V^*)\cong V^*
$$
of (ii). Let $M=\Hom(V,A)$. We identify $V^*$ with a subspace of $M$.
If $e^*_1,\ldots,e^*_n$ and $e_1,\ldots,e_n$ are dual bases for $V^*$ and $V$
then $e^*_1,\ldots,e^*_n$ is also a basis for $M$ over $A$, and
$\si=\sum_j\si(e_j)e_j^*$ for each $\si\in M$. Let $\si_1,\ldots,\si_n$ be a
basis for $D(V)$. By Corollary 1.3 $\si_1,\ldots,\si_n$ is a basis for $M$
over $A$. It follows that the matrix $X$ whose entry in the $i$th row, $j$th
column is $\si_i(e_j)$, $\,1\le i,j\le n$, is invertible in the ring $\Mat_n(A)$.
Given $a_1,\ldots,a_n\in A$ such that $\sum_ia_i\si_i(v)=0$ for all $v\in V$, one
then has necessarily $a_1=\ldots=a_n=0$. In particular, if $\tau_1,\ldots,\tau_n\in D(U)$
and $\sum_i\tau_i(u)\si_i(v)=0$ for all $u\in U$ and $v\in V$, then
$\tau_1=\ldots=\tau_n=0$. This shows that the map in (iii) is injective
(this map is well-defined since the multiplication $A\otimes A\to A$ is a
morphism in $H\mo$). By comparison of dimensions the map in (iii) is bijective.

There is a morphism $k\to V\otimes V^*$ in $H\mo$ given by
$1\mapsto\sum_je_j\otimes e_j^*$. The duality functor produces a map
$D(V)\otimes D(V^*)\cong D(V\otimes V^*)\to D(k)\cong A^H=k$ under which
$\si\otimes\tau\mapsto\sum_j\si(e_j)\tau(e_j^*)$. This shows that
$\<\si,\tau\>\in k$ for all $\si\in D(V)$, $\tau\in D(V^*)$. If
$\<\si_i,\tau\>=0$ for all $i=1,\ldots,n$ then $\tau=0$ by the invertibility
of the matrix $X$ above. Since $\dim D(V)=n=\dim D(V^*)$ according to (ii), the pairing in
(iv) is nondegenerate. If now $\tau_1,\ldots,\tau_n$ is a basis of $D(V^*)$,
dual to $\si_1,\ldots,\si_n$ with respect to this pairing, and $Y$ is the
$n\times n$-matrix containing $\tau_j(e_i^*)$ in the $i,j$-position, then
$XY$ equals the identity matrix. Hence so is $YX$ too, since $X$ is
invertible. This gives $\sum_l\tau_l(e_i^*)\si_l(e_j)=\de_{ij}$,
yielding (v) for $\xi=e_i^*$ and $v=e_j$. Since the expressions in (v) are
linear in $\xi$ and $v$, the equality is fulfilled in general.
\endproof

\Remark.
Both (iv) and (v) are consequences of (iii) in view of the abstract
characterization of dual objects in tensor categories as presented, e.g., in
\cite{Kas, XIV.2}. By the above proposition, for any $U',V'\in E\mo$ there
is a $E$-module structure on $U'\otimes V'$ coming from the isomorphism
$U'\otimes V'\cong D\bigl(D^{-1}(U')\otimes D^{-1}(V')\bigr)$. If also
$W'\in E\mo$, then the canonical map $U'\otimes(V'\otimes W')\to(U'\otimes
V')\otimes W'$ is a morphism in $E\mo$. Thus $E\mo$ becomes a tensor
category in which every object has a left and a right duals. This additional
structure on $E\mo$ corresponds to a Hopf algebra structure on $E$. A more
explicit argument is given below.
\endremark

\proclaim
Lemma 1.5.
Given $\ph\in E$, there are unique elements $\sum_i\ph_i'\otimes\ph_i''\in E\otimes E$
and $\sum_j\ph_j'\otimes\ph_j''\otimes\ph_j'''\in E\otimes E\otimes E$
such that $\ph(ab)=\sum_i\ph_i'(a)\ph_i''(b)$ for all $a,b\in A$ and
$\ph(abc)=\sum_j\ph_j'(a)\ph_j''(b)\ph_j'''(c)$ for all $a,b,c\in A$.
\endproclaim

\Proof.
Note that $E=D(A)$. By Proposition 1.3(iii) $E\otimes E\cong D(A\otimes A)$
and $E\otimes E\otimes E\cong D(A\otimes A\otimes A)$. It remains to notice
that the maps $a\otimes b\mapsto\ph(ab)$ and $a\otimes b\otimes c\mapsto\ph(abc)$
are elements of $D(A\otimes A)$  and $D(A\otimes A\otimes A)$, respectively.
\endproof

\proclaim
Theorem 1.6.
There is a unique comultiplication on $E$ with respect to which $E$ is a Hopf algebra
and $A$ is a left $E$-module algebra. The algebra $A$ is $E^*$-Galois and
the canonical map $H\to\End_EA$ is bijective.
\endproclaim

\Proof.
Define $\De:E\to E\otimes E$ by assigning to each $\ph\in E$ the element of
$E\otimes E$ described in Lemma 1.5. Both $(\De\otimes\id)\De(\ph)$ and
$(\id\otimes\De)\De(\ph)$ correspond to the map
$a\otimes b\otimes c\mapsto\ph(abc)$ in $D(A\otimes A\otimes A)$. Hence
$\De$ is coassociative. Next, we have $\ph(1)\in A^H=k$ for all $\ph\in E$.
Define $\ep:E\to k$ setting $\ep(\ph)=\ph(1)$. Substituting either $a=1$ or
$b=1$ in the identity $\ph(ab)=\sum_{(\ph)}\ph'(a)\ph''(b)$, we see that
$\ep$ is a counit for $\De$. If $\ph,\psi\in E$ then
$$
(\ph\psi)(ab)=\ph\bigl(\textsum_{(\psi)}\psi'(a)\psi''(b)\bigr)
=\textsum_{(\ph),(\psi)}\ph'\psi'(a)\ph''\psi''(b)
$$
for all $a,b\in A$, whence $\De(\ph\psi)=\De(\ph)\De(\psi)$. One has also
$\De(\id_A)=\id_A\otimes\id_A$ so that $\De$ is a homomorphism of unital algebras.

We will check next that the bialgebra $E$ has an antipode. Let
$E^+=\{\ph\in E\mid\ep(\ph)=0\}$ denote the augmentation ideal of $E$.
Since $A$ is a cyclic free $E$-module, $E^+\!A$ is a subspace of codimension 1
in $E$. Take $\chi\in A^*$ such that $\ker\chi=E^+\!A$. Since the actions of
$H$ and $E$ on $A$ commute with each other, $E^+\!A$ is $H$-stable. Then
$A/E^+\!A$ is a onedimensional $H$-module on which $H$ operates via an algebra
homomorphism $\al:H\to k$. So $\chi(ha)=\al(h)\chi(a)$ for all $h\in H$ and
$a\in A$. By Proposition 1.2(ii) the bilinear form $(a,b)=\chi(ab)$ is
nondegenerate on $A$. Hence for each $\ph\in E$ there exists the adjoint
linear transformation $S(\ph):A\to A$ such that $(\ph a,b)=(a,S(\ph)b)$ for
all $a,b\in A$. Let $k_\al=k$ as a vector space, and let $H$ operate in $k_\al$
via $\al$. The map $A\otimes A\to k_\al$, $\,a\otimes b\mapsto\chi(ab)$ is
then a morphism in $H\mo$. By (1.2) it corresponds to an $H$-module map
$\nu:A\to\*\!\!A\otimes k_\al$. Forgetting the $H$-module structure, we may
identify the target of $\nu$ with the dual of $A$, and then $\nu(b)(a)=\chi(ab)$
for $a,b\in A$, so that $\nu$ is bijective. Next,
$S(\ph)=\nu^{-1}\circ(\*\!\ph\otimes\id)\circ\nu$ where $\*\!\ph$
is the $H$-module endomorphism of $\*\!\!A$ defined by $\xi\mapsto\xi\circ\ph$.
It follows that $S(\ph)\in\End_HA$. We have thus constructed a map $S:E\to E$. Now
$$
\ep(\ph)\chi(ab)=\chi\bigl(\ph(ab)\bigr)
=\chi\bigl(\textsum_{(\ph)}\ph'a\cdot\ph''b\bigr)
=\chi\bigl(a\cdot\textsum_{(\ph)}S(\ph')\ph''b\bigr)
$$ for all $a,b$, and the
nondegeneracy of the bilinear form entails
$\sum_{(\ph)}S(\ph')\ph''b=\ep(\ph)b$. Hence $S$ is a left inverse of
$\id_E$ in the convolution algebra $\Hom(E,E)$. Since left invertible
elements in finite dimensional associative algebras are invertible on both
sides, $S$ is the antipode of $E$.

Corollary 1.3 shows for $V=A$ that the map $A\otimes E\to\End A$ defined in (1.3) with
$E$ in place of $H$ is bijective. Hence $A$ is $E^*$-Galois. The bijectivity of
the map $H\to\End_EA$ is a general property of Morita dualities.
\endproof

\Remarks.
The Hopf algebra $L=L(A,H^*)$ constructed by Schauenburg is characterized by
the properties that $A$ is left $L$-Galois and the coaction of $L$ on $A$ commutes
with that of $H^*$ \cite{Scha96, Th.~3.5}. It follows that $L\cong(E\op)^*$. The
existence of antipode in $L$ is a consequence a result in \cite{Scha97} according to
which every bialgebra coacting on a Galois algebra is in fact a Hopf algebra.

There is a well-known characterization of Galois algebras as crossed products
\cite{Doi86, Th.~9, Th.~11}, \cite{Bl89, Th.~1.18}. According to these results one can
identify $A$ with $H^*$ as a left $H$-module, and the multiplication in $A$ is
then obtained by twisting the multiplication in $H^*$ so that
$\xi\mathbin{\cdot_J}\eta=\sum_{(\xi),(\eta)}\<\xi'\otimes\eta',J\>\xi''\eta''$ for
$\xi,\eta\in H^*$ where $J\in H\otimes H$ is an invertible cocycle.
This realization of $A$ provides a particular choice of an algebra isomorphism
$E\cong H\op$ with respect to which the comultiplication in $E$ corresponds to
the $J$-twist $h\mapsto J\De(h)J^{-1}$ of the comultiplication in $H$ as defined by
Drinfeld \cite{Dr89b}. From the viewpoint of Drinfeld's twists some properties of the
transformed Hopf algebra were studied in \cite{Ali01}. The nondegenerate
twists defined in \cite{Ali01} correspond to central simple Galois algebras.
\endremark

\proclaim
Lemma 1.7.
An element $\ph\in E$ is grouplike if and only if $\ph$ is an automorphism
of $A$. Similarly, $h\in H$ is grouplike if and only if $h$ operates on $A$ as
an automorphism.
\endproclaim

This is immediate from Lemma 1.5 and the symmetry between $H$ and $E$.

\proclaim
Lemma 1.8.
Let $V\in H\mo$. The pairings of Proposition {\rm1.4(iv)} induce isomorphisms
$D(V^*)\cong\*\!D(V)$ and $D(\*V)\cong D(V)^*$ in $E\mo$.
\endproclaim

\Proof.
The pairing $D(V)\times D(V^*)\to k$ is $E$-invariant in the sense that
$\ep(\ph)\<\si,\tau\>=\sum_{(\ph)}\<\ph'\si,\ph''\tau\>$ for all
$\si\in D(V)$, $\tau\in D(V^*)$ and $\ph\in E$. A standard transformation of
this identity gives $\<\si,S(\ph)\tau\>=\<\ph\si,\tau\>$, whence the
first isomorphism. Since $(\*V)^*=V$, there is also an $E$-invariant
pairing $D(\*V)\times D(V)\to k$.
\endproof

\section
2. A correspondence between triangular structures

A {\it quasitriangular structure} on a Hopf algebra $H$ is an
element $R=\sum_i\,r'_i\otimes r''_i\in H\otimes H$ satisfying the following conditions:
\mlines
$$
\De\op(h)R=R\De(h)\qquad{\fam0for\ all\ }h\in H,\eqno(2.1)
$$$$
(\ep\otimes\id)(R)=1,\qquad\qquad(\id\otimes\ep)(R)=1,\eqno(2.2)
$$$$
(\De\otimes\id)(R)=R_{13}R_{23},\qquad
(\id\otimes\De)(R)=R_{13}R_{12}\eqno(2.3)
$$
\mlines
where $R_{pq}$ for $1\le p,q\le n$ denote the image of $R$ under the algebra
homomorphism $H\otimes H\to H^{\otimes n}$ which identifies
$H\otimes 1$ and $1\otimes H$ with the $p$th and $q$th tensorands of
$H^{\otimes n}$, respectively. As shown in \cite{Rad92, Lemma 1}, the above
conditions imply that $R$ is invertible (with $R^{-1}=\sum_i\,S(r'_i)\otimes
r''_i$) in accordance with the original definition due to Drinfeld \cite{Dr87}.
If $R^{-1}=R_{21}$, then $R$ is called {\it triangular}.

The map $f_r:H^*\to H$ defined by $f_r(\xi)=\sum_i\,\<\xi,r'_i\>r''_i$ for $\xi\in
H^*$ is a homomorphism of algebras and an antihomomorphism of coalgebras.
Similarly, $f_l:H^*\to H$ defined by $f_l(\xi)=\sum_i\,\<\xi,r''_i\>r'_i$ is
an antihomomorphism of algebras and a homomorphism of coalgebras.
As a consequence $R_l=\Im f_l$ and $R_r=\Im f_r$ are Hopf subalgebras of $H$.
If in the expression $\sum_i\,r'_i\otimes r''_i$ the elements
$\{r'_i\}$ form a basis for $R_l$, then the elements $\{r''_i\}$ form a basis for
$R_r$. Define the {\it rank} of $R$ as $\rk R=\dim R_l=\dim R_r$ \cite{Rad93}. By
Nichols and Zoeller's freeness theorem \cite{Nich89} $\rk R\mid\dim H$ when $H$ is
finite dimensional. We say that $R$ is of {\it maximal rank} if $\rk R=\dim H$.

A quasitriangular structure $R$ on $H$ determines a braiding on the tensor
category $H\mo$ given by the maps
$$
c_{UV}:U\otimes V\iso V\otimes U,\qquad\qquad
u\otimes v\mapsto\textsum_i\,r''_iv\otimes r'_iu
$$
(see \cite{Kas, VIII.3}). This makes $H\mo$ into a braided tensor category as defined
by Joyal and Street \cite{Joy93}.

\proclaim
Lemma 2.1.
Suppose that $R$ is a quasitriangular structure on H and $A$ is a
left $H$-module algebra. Then the formula $T(a,b)=\sum_i(r''_ib)(r'_ia)$ for $a,b\in A$ defines a new
associative multiplication on $A$ such that $T(a,1)=T(1,a)=a$ for all $a\in A$.
The map $T$ is $H$-invariant, so that $A_R=(A,T)$ is a left $H$-module algebra.
If $A$ is $H^*\!$-Galois then so is $A_R$ too and $\End_HA_R\cong(\End_HA)\cop$
as Hopf algebras.
\endproclaim

\Proof.
The map $A\otimes A\to A$ given by $a\otimes b\mapsto T(a,b)$ is the composite of the braiding map
$c_{AA}:A\otimes A\to A\otimes A$ and the original multiplication $A\otimes A\to A$, both of which
are homomorphisms of $H$-modules. The unity property of 1 follows from (2.2). The associativity of
$T$ is a consequence of formal properties of the braiding ($A_R$ and $A$ are opposite algebras in the
braided tensor category of $H$-modules). A straightforward verification is as
follows. If $a,b,c\in A$ then
$$
\eqalign{
T(a,T(b,c))&{}=\textsum_{i,l}\,r''_i\bigl((r''_lc)(r'_lb)\bigr)\cdot(r'_ia)
=\textsum_{i,j,l}\,(r''_jr''_lc)(r''_ir'_lb)(r'_ir'_ja)\,,\cr
T(T(a,b),c)&{}=\textsum_{i,l}\,(r''_ic)\cdot r'_i\bigl((r''_lb)(r'_la)\bigr)
=\textsum_{i,j,l}\,(r''_ir''_jc)(r'_ir''_lb)(r'_jr'_la)\cr
}
$$
by (2.3). The two expressions above are equal as follows from the well known identity
$R_{12}R_{13}R_{23}=R_{23}R_{13}R_{12}$ (see \cite{Dr89a}) explicitly
written as
$$
\textsum_{i,j,l}\,r'_ir'_j\otimes r''_ir'_l\otimes r''_jr''_l
=\textsum_{i,j,l}\,r'_jr'_l\otimes r'_ir''_l\otimes r''_ir''_j\,.
$$

Consider the canonical homomorphism $\pi:A_R\#H\to\End A$. Thus
$\pi(a\#h)(b)=T(a,hb)=\sum_i\,(r''_ihb)(r'_ia)$ for $a,b\in A$ and $h\in H$.
Let $R^{-1}=\sum_j\r'_j\otimes\r''_j$. Then
$$
\textsum_j\,\pi(\r'_ja\#\r''_j)(b)=\textsum_{i,j}\,(r''_i\r''_jb)(r'_i\r'_ja)=ba
$$
since $\sum_{i,j}\,r'_i\r'_j\otimes r''_i\r''_j=1\otimes1$. This shows that $\Im\pi$
contains the subalgebra of right multiplication operators $b\mapsto ba$ in $A$.
On the other hand $\Im\pi$ contains also the subalgebra of operators in $\End A$ given
by actions of elements in $H$. By the bijectivity of map (1.4) these two subalgebras generate the whole
$\End A$ whenever $A$ is $H^*$-Galois. In this case $\pi$ is surjective,
hence bijective by dimension considerations, so that $A_R$ is $H^*$-Galois.
Let $E=\nobreak\End_HA$. Then $\End_HA_R=E$ as algebras since $A_R=A$ as
$H$-modules. Finally, $A_R$ is a left $E\cop$-module algebra as for $\ph\in E$ we have
$$
\ph\bigl(T(a,b)\bigr)=\textsum_{(\ph),i}\,(\ph'r''_ib)(\ph''r'_ia)
=\textsum_{i,(\ph)}\,(r''_i\ph'b)(r'_i\ph''a)=\textsum_{(\ph)}\,T(\ph''a,\ph'b).
\vadjust{\vskip-20pt}
$$
\hfill
\endproof

Further on we make the same assumptions as in section 1, that is,
$A$ is $H^*\!$-Galois and $E=\End_HA$.

\proclaim
Proposition 2.2.
There is a bijective correspondence between the quasitriangular structures on
$H$ and $E$ under which $R=\sum_ir'_i\otimes r''_i$ and
$\R=\sum_m\rho'_m\otimes\rho''_m$ correspond to one another if and only if
$$
\textsum_i\,(r''_ib)(r'_ia)=\textsum_m\,(\rho''_ma)(\rho'_mb)\qquad
{\fam0for\ all\ }a,b\in A.\eqno(2.4)
$$
Furthermore, $\R$ is triangular if and only if so is $R$.
If $c$ and $\c$ denote the bradings of the categories $H\mo$ and $E\mo$
determined by $R$ and $\R$, respectively, then for each pair $U,V\in H\mo$
there is a commutative diagram
$$
\vcenter{\diagram{
D(V\otimes U)&\longmapright8{D(c_{UV})}&D(U\otimes V)\cr
\diagramskip
\mapup{}{}&&\mapup{}{}\cr
\diagramskip
D(V)\otimes D(U)&\longmapright8{\c_{D(V),D(U)}}&D(U)\otimes D(V)\cr
}}\eqno(2.5)
$$
where the vertical arrows are the bijections of Proposition {\rm1.4(iii)}.
\endproclaim

\Proof.
According to Proposition 1.4 $D(A\otimes A)\cong D(A)\otimes D(A)=E\otimes E$. Given
$R$, the map $a\otimes b\mapsto\sum_i\,(r''_ib)(r'_ia)$ is an element of
$D(A\otimes A)$ by Lemma 2.1. It follows that there exists a unique element
$\R\in E\otimes E$ satisfying identity (2.4). We will check that $\R$ is a
quasitriangular structure on $E$. If $\ph\in E$ then, applying $\ph$ to both sides
of (2.4), we deduce
$$
\textsum_{(\ph),i}\,(\ph'r''_ib)(\ph''r'_ia)
=\textsum_{(\ph),m}\,(\ph'\rho''_ma)(\ph''\rho'_mb)\,.
$$
Since the actions of $H$ and $E$ on $A$ commute with each other,
the left hand side can be rewritten as
$$
\textsum_{(\ph),i}\,(\ph'r''_ib)(\ph''r'_ia)
=\textsum_{i,(\ph)}\,(r''_i\ph'b)(r'_i\ph''a)
=\textsum_{m,(\ph)}\,(\rho''_m\ph''a)(\rho'_m\ph'b)\,,
$$
again by (2.4). Thus
$\sum_{(\ph),m}\,\rho''_m\ph''\otimes\rho'_m\ph'
=\sum_{(\ph),m}\,\ph'\rho''_m\otimes\ph''\rho'_m$ as both sides of this equality
correspond to the same element in $D(A\otimes A)$. Swapping the tensorands
yields $\R\De(\ph)=\De\op(\ph)\R$.

Taking $b=1$ in (2.4), we get
$\sum_m\ep(\rho'_m)\rho''_ma=\sum_i\ep(r''_i)r'_ia=a$ for all $a$.
Hence $(\ep\otimes\id)(\R)=1$. The substitution $a=1$ in (2.4)
gives similarly $(\id\otimes\ep)(\R)=1$. For $a,b,c\in A$ we have
$$
\eqalign{
\textsum_m\,(\rho''_ma)\rho'_m(bc)=\textsum_i\,r''_i(bc)(r'_ia)
&{}=\textsum_{i,j}\,(r''_jb)(r''_ic)(r'_ir'_ja)\cr
&{}=\textsum_{j,m}\,(r''_jb)(\rho''_mr'_ja)(\rho'_mc)\cr
&{}=\textsum_{j,m}\,(r''_jb)(r'_j\rho''_ma)(\rho'_mc)\cr
&{}=\textsum_{m,n}\,(\rho''_n\rho''_ma)(\rho'_nb)(\rho'_mc).\cr
}
$$
The map $A\otimes A\otimes A\to A$ given by the rule
$a\otimes b\otimes c\mapsto\sum_m(\rho''_ma)\rho'_m(bc)$ is clearly a
morphism in $H\mo$. By Proposition 1.4
$D(A\otimes A\otimes A)\cong E\otimes E\otimes E$. It follows from the
displayed equality that
$(\id\otimes\De)(\sum_m\rho''_m\otimes\rho'_m)
=\sum_{m,n}\rho''_n\rho''_m\otimes\rho'_n\otimes\rho'_m$, and thus
$(\De\otimes\id)(\R)
=\sum_{m,n}\rho'_n\otimes\rho'_m\otimes\rho''_n\rho''_m
=\R_{13}\R_{23}$. Taking $a$ in (2.4) to be a product of two elements, we
deduce similarly that $(\id\otimes\De)(\R)=\R_{13}\R_{12}$.
We have thus constructed a map $R\mapsto\R$ in one direction. By symmetry
between $H$ and $E$ there is also a map in the opposite direction which is
clearly the inverse one. Next,
$$
\eqalign{
\textsum_{i,j}\,(r''_ir'_ja)(r'_ir''_jb)
=\textsum_{j,m}\,(\rho''_mr''_jb)(\rho'_mr'_ja)
&{}=\textsum_{j,m}\,(r''_j\rho''_mb)(r'_j\rho'_ma)\cr
&{}=\textsum_{m,n}\,(\rho''_n\rho'_ma)(\rho'_n\rho''_mb).\cr
}
$$
If $R_{21}R=1\otimes1$ then the left hand side is equal to $ab$ for all $a,b$,
whence $\R_{21}\R=\sum_{m,n}\,\rho''_n\rho'_m\otimes\rho'_n\rho''_m=1\otimes1$.

Let $\tau\in D(U)$, $\si\in D(V)$, $u\in U$, $v\in V$. The commutativity of
(2.5) is verified by the following calculation:
$$
\textsum_i\,\si(r''_iv)\tau(r'_iu)
=\textsum_i\,r''_i\si(v)r'_i\tau(u)
=\textsum_m\,\rho''_m\tau(u)\rho'_m\si(v)\,.
\vadjust{\vskip-20pt}
$$
\hfill
\endproof

Until the end of this section we assume that $R$ and $\R$ are mutually
corresponding quasitriangular structures on $H$ and $E$. The next result is a
special case of Theorem 2.5. However it is much easier to handle this special
case which will be eventually of main interest to us.

\proclaim
Proposition 2.3.
Suppose that $A$ is central simple. If $H$ is cocommutative and $R=1\otimes1$
then $\rk\R=\dim A$.
\endproclaim

\Proof.
Under the stated hypotheses equation (2.4) can be rewritten as
$$
ba=\textsum_m\,(\rho''_ma)(\rho'_mb).
$$
It shows that the right multiplication by $a$ in $A$ is given by action of the
element $\sum_m\,\rho''_ma\#\rho'_m\in A\#E$. Since
$\sum_m\,\rho'_m\otimes\rho''_m\in\R_l\otimes\R_r$, we have
$$
\textsum_m\,\rho''_ma\#\rho'_m\in A\#\R_l.
$$
Then the image of the canonical homomorphism $A\otimes A\op\to\End A$ is
contained in the image of $A\#\R_l$. Since $A$ is central simple, one has
$A\otimes A\op\cong\End A$, whence the map $A\#\R_l\to\End A$ is surjective.
The latter map is the restriction of the canonical homomorphism $A\#E\to\End A$
which is bijective since $A$ is $E^*$-Galois. It follows that $\R_l=E$, and so
$\rk\R=\dim\R_l=\dim E=\dim A$.
\endproof

\proclaim
Lemma 2.4.
There is a homomorphism of algebras $\Phi:A_R\otimes A_R\op\to A\#(H\otimes
E)$ such that
$$
a\otimes1\mapsto u_a=\textsum_m\,\rho''_ma\#(1\otimes\rho'_m),\qquad\quad
1\otimes a\mapsto v_a=\textsum_i\,r''_ia\#(r'_i\otimes1)
$$
for $a\in A$. If the algebra $A_R$ is central simple then\/
$\dim A\mid(\rk R)(\rk\R)$.
\endproclaim

\Proof.
We may regard $A$ as an $(H\otimes E)$-module algebra since the actions of $H$
and $E$ on $A$ commute. Let $\pi:A\#(H\otimes E)\to\End A$ be the
canonical homomorphism. If $a,b\in A$ then $\pi(v_a)(b)=T(b,a)$ where $T$ is the
multiplication in $A_R$. By (2.4) we have
also $\pi(u_a)(b)=T(a,b)$. In other words, $\pi(v_a)$ and $\pi(u_a)$
coincide respectively with the right and the left multiplication operators
in $A_R$. Identify $H$ and $E$ with the subalgebras $H\otimes 1$ and
$1\otimes E$ of $H\otimes E$. The restriction of $\pi$ to each of these two
is an isomorphism since $A$ is both $H^*\!$-Galois and $E^*\!$-Galois. Note that
$u_a\in A\#E$ and $v_a\in A\#H$. It follows that $u_au_b=u_{T(a,b)}$ and
$v_av_b=v_{T(b,a)}$ for all $a,b\in A$. Next,
$$
\eqalign{
v_au_b&{}=\,\textsum_{i,j,n}\,\,(r''_ir''_ja)(r'_i\rho''_nb)\#(r'_j\otimes\rho'_n),\cr
u_bv_a&{}=\textsum_{j,m,n}\,(\rho''_m\rho''_nb)(\rho'_mr''_ja)\#(r'_j\otimes\rho'_n).\cr
}
$$
In view of (2.4) $v_au_b=u_bv_a$. Hence the map $\Phi$ defined by the rule
$a\otimes b\mapsto u_av_b$ is a homomorphism of algebras.

We see that $\Im\Phi\subset B=A\#(R_l\otimes\R_l)$. Suppose that $A_R$ is
central simple. Then so are $A_R\op$ and $A_R\otimes A_R\op$ as well. Hence
$\Phi$ is an embedding and $B\cong\Im\Phi\otimes C$ where $C$ is the
centralizer of $\Im\Phi$ in $B$ by \cite{Jac, Th.~4.7}. It follows that $\dim\Im\Phi\mid\dim
B$. As $\dim\Im\Phi=(\dim A)^2$ and $\dim B=(\dim A)(\rk R)(\rk\R)$, the
final claim of the lemma is proved.
\endproof

\proclaim
Theorem 2.5.
If $A$ is central simple then\/ $\dim A\mid(\rk R)(\rk\R)$.
\endproclaim

\Proof.
The passage from $A$ to $A_R$ can be reversed with the same construction.
As is easily checked, $R_{21}^{-1}$ is another quasitriangular structure on
$H$ \cite{Maj, Exercise 2.1.3}. Letting $R^{-1}=\sum_j\r'_j\otimes\r''_j$, the multiplication in the
algebra $(A_R)_{R_{21}^{-1}}$ is given by the formula
$$
(a,b)\mapsto\textsum_j\,T(\r'_jb,\r''_ja)
=\textsum_{i,j}\,(r''_i\r''_ja)(r'_i\r'_jb)=ab,
$$
which is the original multiplication in $A$. By Lemma 2.1 $A_R$ is $H^*\!$-Galois and
$\End_HA_R\cong E\cop$. Let $\R^{-1}=\sum_n\rhot'_n\otimes\rhot''_n$. As
$$
ab=\textsum_{m,n}\,(\rho''_m\rhot''_na)(\rho'_m\rhot'_nb)
=\textsum_n\,T(\rhot''_na,\rhot'_nb)
$$
for all $a,b\in A$, formula (2.4) applied in case of the $H^*\!$-Galois
algebra $A_R$ shows that $\R^{-1}$ is the quasitriangular structure on
$E\cop$  corresponding to $R_{21}^{-1}$. Since
$R^{-1}=\sum_i\,S(r'_i)\otimes r''_i$ by \cite{Dr89a}, we have
$\rk R_{21}^{-1}=\rk R$. Similarly, $\rk\R^{-1}=\rk \R$. It remains to apply Lemma
2.4 with $A_R$ in place of $A$ and $R_{21}^{-1}$ in place of $R$.
\endproof

\Remark.
At the other extreme when $A$ is commutative (and $H$ is necessarily
cocommutative), one has $\rk R=\rk\R$.
\endremark

\section
3. The Miyashita-Ulbrich action and the square of the antipode

Given an $H$-Galois algebra, there is a canonically defined right action of
$H$ on the algebra under consideration \cite{Ul82}, \cite{Doi89}.
The next proposition restates several facts from \cite{Doi89, Th.~3.4},
with $H$ and $H^*$ interchanged by our conventions. So let $A$ be again an $H^*$-Galois algebra.

\proclaim
Proposition 3.1.
Suppose $B$ is an algebra containing $A$ as a subalgebra.
There is a unique right $H^*\!$-module structure on $B$ such that
$$
ba=\textsum_i\,(h_ia)(bh_i^*)\qquad{\fam0for\ all\ }a\in A
{\fam0\ and\ }b\in B\eqno(3.1)
$$
where $\{h_i\}$ and $\{h_i^*\}$ are dual bases of $H$ and $H^*$, respectively.
It makes $B$ into an $H^*\!$-module algebra, and the subalgebra of invariants
$B^{H^*}$ coincides with the centralizer of $A$ in $B$.
Let $\xi\in H^*$. If $\sum_jv_j\otimes w_j\in A\otimes A$ is the preimage of $1\otimes\xi$ under
bijection {\rm(1.5)}, so that $\sum_jv_j(hw_j)=\<\xi,h\>1$ for all $h\in H$,
then
$$
\vadjust{\vskip-8pt}
b\xi=\textsum_j\,v_jbw_j\qquad{\fam0for\ all\ }b\in B.\eqno(3.2)
$$
\endproclaim

In particular, $A$ is a right $H^*\!$-module algebra in a natural way.
The $H^*\!$-comodule structure on $A$ is given by the map
$a\mapsto\sum\,h_ia\otimes h_i^*$. Hence (3.1) is exactly the exchange rule
from \cite{Ul82}, \cite{Doi89}. The symmetry between $H$ and $H^*$ is apparent in
(3.1). Formula (3.2) was used by Ulbrich \cite{Ul82} to define the module
structure explicitly. This formula reflects the fact that the preimage of
$1\otimes H^*$ under bijection (1.5) is a subalgebra of the tensor product algebra $A\otimes A\op$
isomorphic to $(H^*)\op$. This can be further refined as follows:

\proclaim
Lemma 3.2.
Bijection $\ga$ given by {\rm(1.5)} is an isomorphism between two algebra structures:
$A\otimes A\op\iso A\#(H^*)\op$. If $B$ is an algebra containing $A$ as a
subalgebra and $\pit:A\#(H^*)\op\to\End B$ is defined by
$\pit(a\#\xi)b=a(b\xi)$, then
$$
(\pit\circ\ga)(c\otimes a)(b)=cba,\qquad
{\fam0where\ }a,c\in A,\ b\in B,\ \xi\in H^*.
\vadjust{\vskip-4pt}
$$
\endproclaim

\Proof.
As $\ga(c\otimes a)=\sum_ic(h_ia)\#h_i^*$, the formula for $\pit\circ\ga$
follows at once from (3.1). We see that $\pit\circ\ga$ is an algebra
homomorphism. Since $\pit$ is also an algebra homomorphism, we can conclude
that so is $\ga$ too provided that $\pit$ is injective. Since $\ga$ is a
bijection, $\pit$ is injective if and only if so is $\pit\circ\ga$. Thus it
suffices to embed $A$ into a central simple algebra $B$. For instance, we may take $B=\End A$.
\endproof

\proclaim
Proposition 3.3.
The algebra $A$ is separable if and only if $A$ is $H^*$-projective. The
following conditions are equivalent:

\item(i)
$A$ is central simple,

\item(ii)
$A$ is left $H$-Galois with respect to the comodule structure
$a\mapsto\sum_i\,h_i\otimes ah_i^*$,

\item(iii)
$A$ is $H^*$-projective and $A^{H^*}=k$.

\endproclaim

This is again just a special case of results in \cite{Ul82}, \cite{Doi89}.
A direct argument based on Lemma 3.2 runs as follows. First note that ``$A$
is $H^*$-projective" implies that $A$ is $D$-projective where we put
$D=A\#(H^*)\op$. In fact $A$ is a $D$-module direct summand of
$M=D\otimes_{(H^*)\op}A$ as the $D$-module map $A\to M$,
$\,a\mapsto(a\#1)\otimes1$, has a left inverse $M\to A$, $\,d\otimes a\mapsto
da$. Clearly $M$ is $D$-projective when $A$ is $H^*$-projective. Conversely,
``$A$ is $D$-projective" implies ``$A$ is $H^*$-projective" since $D$ is a
free left module over $(H^*)\op$. In view of Lemma 3.2 $A$ is
$H^*$-projective if and only if $A$ is $A\otimes A\op$-projective. The
latter condition is the defining property of separable algebras. Next, $A$
is $H$-Galois if and only if $\pit:A\#(H^*)\op\to\End A$ is bijective, and
$A$ is central simple if and only if so is $\pit\circ\ga$. Hence
(i)$\Leftrightarrow$(ii). As stated in Proposition 1.2 condition (ii) implies that $A$ is a
cyclic free $H^*$-module, whence (ii)$\Rightarrow$(iii). If (iii) holds then
$A$ is separable with center $A^{H^*}=k$, which gives (i).

Let $x\in H$ be a {\it left integral}, that is, a nonzero element such that
$hx=\ep(h)x$ for all $h\in H$. Since the space of left integrals is
onedimensional, there is an algebra homomorphism $\al:H\to k$ such that
$xh=\al(h)x$ for all $h\in H$. This defines a grouplike element $\al\in H^*$
with inverse $\al^{-1}=\al\circ S$. We call $\al$ the {\it modular function\/}
on $H$. The Hopf algebra $H$ is {\it unimodular\/} if $\al=\ep$. We will use
left and right actions of $H^*$ on $H$ given by
$$
\xi\rightharpoonup h=\textsum\,\<\xi,h''\>h',\qquad\qquad
h\leftharpoonup\xi=\textsum\,\<\xi,h'\>h''.\eqno(3.3)
$$
Let $E=\End_HA$ as in the preceding sections. Likewise we may regard $H$ as a
subalgebra of $\End A$.

\proclaim
Proposition 3.4.
The linear function $\chi\in A^*$ defined by $\chi(a)=xa$ satisfies
$$
\chi(\ph a)=\ep(\ph)\chi(a),\qquad\qquad
\chi(ha)=\al(h)\chi(a)\eqno(3.4)
$$
for all $\ph\in E$, $h\in H$ and $a\in A$. There exists an automorphism
of finite order $\th:A\to A$ such that
\vadjust{\vskip-6pt}
$$
\vadjust{\vskip-6pt}
\chi(ba)=\chi(a\cdot\th b)\qquad{\fam0\ for\ all\ }a,b\in A.\eqno(3.5)
$$
One has
$$
\th\ph\th^{-1}=S^2(\ph),\qquad\qquad
\th h\th^{-1}=S^2(\al^{-1}\rightharpoonup h\leftharpoonup\al).\eqno(3.6)
$$
\endproclaim

\Proof.
Clearly $xA\subset A^H=k$, so that $\chi$ is well defined. Now $\chi(\ph
a)=\ph(xa)=\chi(a)\ph(1)$, which gives the first formula in (3.4). The second
one is immediate from the definition of $\al$. By Proposition 1.2 the bilinear
form $(a,b)\mapsto\chi(ab)$ on $A$ is nondegenerate. The existence of an
automorphism $\th$ satisfying (3.5) is a general property of Frobenius
algebras ($\th$ is the Nakayama automorphism of $A$). Since
$(hb)a=\sum_{(h)}h'\bigl(b\cdot S(h'')a\bigr)$, it follows from (3.4) that
$$
\chi\bigl((hb)a\bigr)=\chi\bigl(b\cdot S(h\leftharpoonup\al)a\bigr)\,.
\eqno(3.7)
$$
Applying (3.5) to both sides, this can be rewritten as
$$
\chi\bigl(a\cdot\th(hb)\bigr)
=\chi\bigl(S(h\leftharpoonup\al)a\cdot\th b\bigr)\,.
$$
Observe that $S(h\leftharpoonup\al)\leftharpoonup\al=S(\al^{-1}\rightharpoonup
h\leftharpoonup\al)$. Using (3.7), the right hand side of the displayed
equality above can be rewritten as
$\chi\bigl(a\cdot S^2(\al^{-1}\rightharpoonup h\leftharpoonup\al)\th b\bigr)$.
The nondegeneracy of the bilinear form $(a,b)\mapsto\chi(ab)$ entails
$$
\th(hb)=S^2(\al^{-1}\rightharpoonup h\leftharpoonup\al)\th b.
$$
A similar computation with $\ep$ in place of $\al$ proves the first formula in (3.6).
By \cite{Rad76} the antipode of the Hopf algebra $E$ has finite order, say
$e$. It follows then from (3.6) that $\th^e$ centralizers all $\ph\in E$,
that is, $\th^e\in\End_EA=H$. Since $\th^e$ is an automorphism of $A$, it is a
grouplike element in $H$ by Lemma 1.7. Since the grouplike elements in $H$ constitute a
finite group, $\th^e$ has finite order. Then so does $\th$ too.
\endproof

\proclaim
Proposition 3.5.
Suppose that $R=\sum_ir'_i\otimes r''_i$ and $\R=\sum_m\rho'_m\otimes\rho''_m$
are mutually corresponding quasitriangular structures on Hopf algebras $H$ and $E$,
respectively. Denote by $u=\sum_iS(r''_i)r'_i$ and
$\de=\sum_mS(\rho''_m)\rho'_m$ Drinfeld's elements in $H$ and $E$.
Then $\th=\de g_\al^{-1}S(u)^{-1}$ where $g_\al=\sum_i\al(r'_i)r''_i$ is a grouplike element in $H$.
In particular, $\th=\de\in E$ when $H$ is cocommutative and $R=1\otimes1$.
\endproclaim

\Proof.
Apply $\chi$ to both sides of (2.4). As $\chi(\ph a\cdot b)=\chi(a\cdot
S(\ph)b)$ in view of (3.4), we have
$$
\chi(\textsum_m\rho''_ma\cdot\rho'_mb)
=\chi(a\cdot\de b).\eqno(3.8)
$$
Identity $b(ha)=\sum_{(h)}h''\bigl(S^{-1}(h')b\cdot a\big)$ together with (2.3)
give
$$
\textsum_i\,(r''_ib)(r'_ia)
=\textsum_{i,j}\,r'_j\bigl(S^{-1}(r'_i)r''_ir''_jb\cdot a\bigr).
\vadjust{\vskip-4pt}
$$
Here $\sum_iS^{-1}(r'_i)r''_i=S^{-1}(u)=S(u)$ by \cite{Dr89a}. Using (3.4) and then (3.5), we get
$$
\vadjust{\vskip-4pt}
\chi(\textsum_ir''_ib\cdot r'_ia)
=\chi(S(u)g_\al b\cdot a)=\chi(a\cdot\th S(u)g_\al b).\eqno(3.9)
$$
Comparison of (3.8) with (3.9) gives $\de=\th S(u)g_\al$ in $\End A$.
\endproof

\Remark.
Since $S(u)g_\al\in H$, equality (3.6) shows that $S^2(\ph)=\de\ph\de^{-1}$ for all
$\ph\in E$. This is known to hold in any quasitriangular Hopf algebra
\cite{Dr89a}. If $H$ is cocommutative and $R=1\otimes1$ then $\de$ is an
automorphism of $A$, hence a grouplike element in $E$. This is a general
property of triangular Hopf algebras.
\endremark

\proclaim
Lemma 3.6.
Suppose that $H$ is cocommutative and $A$ is central simple. If $\eta\in
A^*$ is a linear function such that $\ker\eta=[A,A]$, then
$\eta(ha)=\ep(h)\eta(a)$ for all $h\in H$ and $a\in A$.
\endproclaim

\Proof.
The action of $H$ on $A$ is {\it inner}, that is, there exists a convolution
invertible linear map $u:H\to A$ such that $ha=\sum_{(h)}u(h')au^{-1}(h'')$.
This is proved in \cite{Bl89, Th.~2.15}, for Hopf Galois extensions; more general
results in \cite{Bea}, \cite{Kop}, \cite{Mas} show that any measuring action of a
coalgebra on a central simple algebra is inner. Using the identity
$\eta(ab)=\eta(ba)$ and the cocommutativity of $H$, we obtain
$$
\eta(ha)=\eta\bigl(a\textsum_{(h)}u^{-1}(h'')u(h')\bigr)=\ep(h)\eta(a).
\vadjust{\vskip-20pt}
$$
\hfill
\endproof

\proclaim
Proposition 3.7.
If $H$ is cocommutative and $A$ is central simple then $\th$ is given by the
Miyashita-Ulbrich action of the modular function $\al\in H^*$.
\endproclaim

\Proof.
Since $A$ is a cyclic free $H$-module and $H$ is a Frobenius algebra, every
simple $H$-module can be embedded in $A$. In particular, $A$ contains a
nonzero element $w$ such that $hw=\al(h)w$ for all $h\in H$. As $Aw$ is a
nonzero submodule of the simple $A\#H$-module $A$, we must have $Aw=A$. This
shows that $w$ is left invertible, hence invertible in $A$. Similarly, $A$
contains an invertible element $v$ with the property that $hv=\al^{-1}(h)v$ for all
$h\in H$. Now $h(vw)=\ep(h)vw$, that is, $vw\in A^H=k$. We may assume after
rescaling that $w=v^{-1}$. Since $v(hw)=\al(h)1$ for all $h$, formula (3.2)
gives $b\al=vbw$ for $b\in A$.

Take $\eta$ as in Lemma 3.6 (note that $[A,A]$ has codimension 1 in $A$). Define $\chi'\in A^*$ by
$\chi'(a)=\eta(av)$. Then $\chi'(ha)=\eta(a\cdot S(h)v)=\al(h)\eta(av)=\al(h)\chi'(a)$ for all $h\in H$ and
$a\in A$. Again by the properties of the $H$-module $A$ the largest
$H$-semisimple factor module of $A$ contains all onedimensional $H$-modules
with multiplicity 1. Hence $\chi'$ is a scalar multiple of $\chi$. Rescaling
$\eta$, we may achieve $\chi'=\chi$. Now
$$
\chi(ba)=\eta(bav)=\eta(avb)=\chi(avbw)=\chi(a\cdot b\al).
$$
The nondegeneracy of the bilinear form $(a,b)\mapsto\chi(ab)$ entails
$\th b=b\al$ for all $b\in A$.
\endproof

\proclaim
Theorem 3.8.
Suppose that $H$ is cocommutative and $A$ is central simple. Let $n$ be the
order of $\al$ in the group of invertible elements of $H^*$. The antipode
$S$ of the Hopf algebra $E=\End_HA$ satisfies $S^{2n}=\id$. If\/ $1$ is the
only central grouplike element in $H$, then the order of $S^2$ equals $n$
exactly.
\endproclaim

\Proof.
By Proposition 3.3 $A$ is $H$-Galois. In particular, $H^*$ operates faithfully in
$A$. This shows together with Proposition 3.7 that $\th$ has order $n$. Then
$S^{2n}=\id$ by (3.6). Suppose that $S^{2e}=\id$ for some $e>0$. Then $\th^e$
is a grouplike element of $H$, as we have seen in the proof of Proposition 3.4.
On the other hand, $\th\in E$ according to Proposition 3.5. Then $\th^e\in H\cap E$,
that is, $\th^e$ belongs to the center of $H$. Under the assumption that $H$
contains no nontrivial central grouplike elements we obtain $\th^e=1$, that is, $n$
divides $e$.
\endproof

\proclaim
Proposition 3.9.
Suppose that $H$ is cocommutative and $A$ is central simple. Let the
triangular structure $\R=\sum_m\,\rho'_m\otimes\rho''_m$ on $E$
correspond to $R=1\otimes1$. Then the Miyashita-Ulbrich action of $E^*$ on $A$ is given by the rule
$$
a\xi=\textsum_m\,\<\xi,\rho''_m\>\rho'_ma\qquad{\fam0for\ }
a\in A{\fam0\ and\ }\xi\in E^*.\eqno(3.10)
$$
\endproclaim

\Proof.
Since the map $f:E^*\to E$ under which $\xi\mapsto\sum_m\<\xi,\rho''_m\>\rho'_m$ is an antihomomorphism
of algebras, formula (3.10) does define a right $E^*$-module structure. By
Proposition 2.3 $\rk\R=\dim E$. We may assume therefore that $\{\rho''_m\}$ is a
basis for $E$. Let $\{\xi_m\}$ be the dual basis for $E^*$. As
$f(\xi_m)=\rho'_m$ for each $m$, formula (2.4) can be rewritten as
$ba=\sum_m\,(\rho''_ma)(b\xi_m)$. This coincides with (3.1) written for $E$ in
place of $H$.
\endproof

\section
4. Generating subspaces

Denote by $H\alg$ (respectively, $H\coalg$) the category of finite dimensional
left $H$-module algebras (respectively, left $H$-module coalgebras). If
$\mu:B\otimes B\to B$ is the multiplication in $B\in H\alg$ then, by passing to duals,
we obtain an $H$-module homomorphism $\mu^*:B^*\to B^*\otimes B^*$ which makes $B^*$ into
a coalgebra. Because of the isomorphisms $(U\otimes V)^*\cong V^*\otimes U^*$
in $H\mo$, we must identify here $B^*\otimes B^*$ with the dual of $B\otimes
B$ by means of the pairing $\<\xi\otimes\eta,\,a\otimes b\>=\<\xi,b\>\<\eta,a\>$
for $\xi,\eta\in B^*$ and $a,b\in B$. Thus the coalgebra constructed is the dual
of $B\op$. One can use right duals quite similarly. In this way we get two functors
$L,R:H\alg\to H\coalg$ for which $B\mapsto(B\op)^*$ and $B\mapsto\*(B\op)$,
respectively. Both of them are equivalences with inverse functors
$C\mapsto(\*C)\op$ and $C\mapsto(C^*)\op$. All parallel constructions apply
in case of the Hopf algebra $E$. In view of Proposition 1.4(iii) the duality $D$
takes $H\alg$ to $E\coalg$ and $H\coalg$ to $E\alg$. If $C\in H\coalg$ then
$D(C)=\Hom_H(C,A)$ is just a subalgebra of the convolution algebra $\Hom(C,A)$.
If $B\in H\alg$ then the comultiplication $\psi\mapsto\sum_{(\psi)}\psi'\otimes\psi''$
in $D(B)=\Hom_H(B,A)$ is characterized by means of the identity
$\psi(bd)=\sum_{(\psi)}\psi'(b)\psi''(d)$ where $b,d\in B$.

Concerning the notations we use here one should keep in mind what part of
structure is affected by them. For instance, $B\op$ indicates the change of
multiplication. As an $H$-module $B\op=B$. In the same way the notations $B^*$
and $\*\!B$ are used to distinguish two module structures without any relevance
to the coalgebra structure.

\proclaim
Proposition 4.1.
Let $B\in H\alg$.

\item(i)
The $H$-invariant elements of the algebras $B\op\otimes A$ and $A\otimes B\op$
constitute their subalgebras isomorphic to $B'=D\bigl((B\op)^*\bigr)$ and
$B''=D\bigl(\*(B\op)\bigr)$, respectively.

\item(ii)
The $E$-module structures on $B'$ and $B''$ correspond to the action of $E$ on
the second tensorand of $B\op\otimes A$ and the first tensorand of $A\otimes B\op$.

\item(iii)
The pairings $D(B)\times B'\to k$ and $B''\times D(B)\to k$ of
Proposition {\rm1.4(iv)} are given by
$$
\<\psi,\si\>=\textsum_i\psi(b_i)a_i,\qquad\qquad
\<\tau,\psi\>=\textsum_ja'_j\psi(b'_j)\eqno(4.1)
$$
where $\si\in B'$ corresponds to $\sum_ib_i\otimes a_i\in(B\op\otimes A)^H$,
$\,\tau\in B''$ corresponds to $\sum_ja'_j\otimes b'_j\in(A\otimes B\op)^H$, and
$\psi\in D(B)$.

\item(iv)
If the algebra $B$ is generated by an $H$-invariant subspace $V$, then
$(B\op\otimes A)^H$ is generated by $(V\otimes A)^H$ and $(A\otimes B\op)^H$ is
generated by $(A\otimes V)^H$.

\endproclaim

\Proof.
By (1.1) $D(V^*)=\Hom_H(V^*,A)\cong(V\otimes A)^H$ for any $V\in H\mo$ since
$\*(V^*)=V$. Similarly, $D(\*V)\cong(A\otimes V)^H$. This applies to $V=B$.
On the other hand, both $B'$ and $B''$ are subalgebras of the convolution
algebra $\Hom(C,A)$, where $C=(B\op)^*$, isomorphic to $B\op\otimes A$.
This shows that $B'$ is embedded into $B\op\otimes A$ (resp., $B''$ into
$A\otimes B\op$) as a subalgebra. If $\ph\in E$ and $\si\in B'$, then
$\ph\si$ is given by the juxtaposition of maps. When interpreting $\si:B^*\to
A$ as an element of $B\otimes A$, this exactly corresponds to the action of
$\ph$ on the second tensorand. Similarly for $B''$. Next we compute the
pairings in (iii). By the hypotheses $\si(\xi)=\sum_i\<\xi,b_i\>a_i$ for
$\xi\in B^*$. Let $\{e_l\}$ and $\{e_l^*\}$ be dual bases of $B$ and
$B^*$. Then
$$
\<\psi,\si\>
=\textsum_l\psi(e_l)\si(e_l^*)
=\textsum_{i,l}\<e_l^*,b_i\>\psi(e_l)a_i
=\textsum_i\psi(b_i)a_i.
\vadjust{\vskip-6pt}
$$
The expression $\<\tau,\psi\>=\sum_l\tau(e_l^*)\psi(e_l)$ is transformed
similarly with a use of equality $\tau(\xi)=\sum_j\<\xi,b'_j\>a'_j$.
Under the hypotheses of (iv) let $G\subset B'$ be the subalgebra
generated by the $E$-submodule $D(V^*)\subset B'$. Then $G$ is stable
under $E$, and so $G\in E\alg$. It follows that
$G'=\bigl(\*\!D^{-1}(G)\bigr)\op$ is a subalgebra of $B$ containing $V$. Since $V$ generates
$B$, we get $G'=B$. Hence $G=B'$. In view of the identifications in (i)
this gives exactly the first statement in (iv). The other one is similar.
\endproof

\proclaim
Proposition 4.2.
The diagrams below are commutative up to natural isomorphisms:
$$
\diagram{
H\alg&\mapright L&H\coalg&\qquad\qquad&H\alg&\mapright R&H\coalg\cr
\mapdown D{}&&\mapdown{}D&&\mapdown D{}&&\mapdown{}D\cr
E\coalg&\mapleft L&E\alg\,,\!&&E\coalg&\mapleft R&E\alg\,.\!\cr
}
$$
\endproclaim

\vskip-6pt
\Proof.
Let $B\in H\alg$ and $C=(B\op)^*$. The pairing $D(B)\times D(C)\to k$
identifies $D(B)$ with the left dual of $D(C)$ in $E\mo$ by Lemma 1.8. We
have to check that the coalgebra structure on $D(B)$ correctly corresponds to the
algebra structure on $D(C)$ under this pairing. Let $\psi\in D(B)$, and let
$\si,\tau\in D(C)$ correspond to $\sum_ib_i\otimes a_i$ and $\sum_jd_j\otimes c_j$
in $(B\op\otimes A)^H$ under identifications of Proposition 4.1. Then
$\sum_{i,j}d_jb_i\otimes a_ic_j$ corresponds to $\si\tau$. By (4.1)
$$
\<\psi,\si\tau\>=\textsum_{i,j}\psi(d_jb_i)a_ic_j
=\textsum_{i,j,(\psi)}\psi'(d_j)\psi''(b_i)a_ic_j
=\textsum_{(\psi)}\<\psi',\tau\>\<\psi'',\si\>,
\vadjust{\vskip-6pt}
$$
that is, $D(B)$ is the dual of $D(C)\op$. The commutativity of the second diagram is verified similarly.
\endproof

Next we will evaluate the composite functors $H\alg\to E\alg$ at certain
$H$-module algebras of particular interest. Note that $H$ is a left $H$-module
algebra with respect to the adjoint action of $H$, and so is $H^*$
with respect to the action $\<h\rightharpoonup\xi,g\>=\<\xi,gh\>$ where
$g,h\in H$ and $\xi\in H^*$. We use notations $H_{\ad}$ and
$H^*_\rightharpoonup$ to specify these particular module structures.

\proclaim
Proposition 4.3.
One has the following isomorphisms in $E\alg${\rm:}

\item(i)
$D\bigl((A\op)^*\bigr)\cong
D\bigl(\*(A\op)\bigr)\cong E^*_\rightharpoonup,$

\item(ii)
if $B=H^*_\rightharpoonup$ then
$D\bigl((B\op)^*\bigr)\cong D\bigl(\*(B\op)\bigr)\cong A,$

\item(iii)
if $B=H_{\ad}$ then
$D\bigl((B\op)^*\bigr)\cong D\bigl(\*(B\op)\bigr)\cong E_{\ad}$.

\endproclaim

\Proof.
According to Proposition 4.2 $D\bigl(\*(A\op)\bigr)\cong\bigl(D(A)^*\bigr)\op$.
Next, $D(A)=E$. The action of $E$ on $D(A)$ is that by left multiplications,
and the coalgebra structure on $E$ defined in section 1 is a particular case
of those structures on $D(B)$ for $B\in H\alg$. The action of $E$ in $E^*$
used here is given by $\<\ph\rightharpoondown\xi,\psi\>=\<\xi,S(\ph)\psi\>$
where $\ph,\psi\in E$ and $\xi\in E^*$. Since $S$ is a Hopf algebra antiautomorphism,
the assignment $\xi\mapsto\xi\circ S$ defines an isomorphism of algebras
$(E^*)\op\to E^*$. Under this bijection
$(\ph\rightharpoondown\xi)\circ S=\ph\rightharpoonup(\xi\circ S)$
as both linear functions produce $\<\xi,S(\ph)S(\psi)\>$ when
evaluated at $\psi$. Thus $(E^*)\op\cong E^*_\rightharpoonup$ in $E\alg$.
To prove the other isomorphism in (i) one should just replace $S$ with
$S^{-1}$ everywhere above.

Replacing $E$ with $H$ in the argument above, we get $(H^*)\op\cong H^*_\rightharpoonup$
in $H\alg$. Hence $\*(B\op)\cong H$ in $H\coalg$ for $B=H^*_\rightharpoonup$. Similarly
$(B\op)^*\cong H$. To prove (ii) it remains to check that $D(H)\cong A$ in $E\alg$. The
assignment $a\mapsto\si_a$ where $\si_a(h)=ha$ for $a\in A$ and $h\in H$ gives
a linear bijection $A\to D(H)$. If $b\in A$ then the equality
$h(ab)=\sum_{(h)}(h'a)(h''b)$ shows that $\si_{ab}=\si_a\si_b$.
The identity $\ph(ha)=h(\ph a)$ for $\ph\in E$ gives $\ph\si_a=\si_{\ph
a}$. Thus the multiplications and the $E$-module structures do correspond
correctly.

Let us prove (iii). We regard $\End A$ as an $H$-module algebra with respect
to the adjoint action $h\triangleright\ph=\sum_{(h)}h'\circ\ph\circ S(h'')$ where
$h\in H$ and $\ph\in\End A$. Note that $E=(\End A)^H$. The embedding
$H_{\ad}\to\End A$ and the map $A\to\End A$ given by left multiplication
operators $a_L$ are morphisms in $H\mo$. Hence so too is the map
$\pi:A\otimes H_{\ad}\to\End A$ defined in (1.3). It follows that $\pi$ maps
$(A\otimes H_{\ad})^H$ bijectively onto $E$. Let $x=\sum_ia_i\otimes h_i$ and
$y=\sum_jc_j\otimes g_j$ be two elements in $(A\otimes H_{\ad})^H$. We get
\vadjust{\vskip-6pt}
$$
\pi(x)\pi(y)=\textsum_i(a_i)_L\circ h_i\circ\pi(y)
=\textsum_i(a_i)_L\circ\pi(y)\circ h_i=\textsum_{i,j}(a_ic_j)_L\circ g_jh_i.
\vadjust{\vskip-6pt}
$$
This shows that $(A\otimes H_{\ad})^H$ regarded as a subalgebra of $A\otimes
H\op$ is isomorphic to $E$ as an algebra under $\pi$. Proposition 4.1(i) with
$B=H_{\ad}$ yields now $D\bigl(\*(B\op)\bigr)\cong E$. If $x$ is as before
and $\ph\in E$ then
$$
\textsum_i(\ph a_i)_L\circ h_i
=\textsum_{i,(\ph)}\ph'\circ(a_i)_L\circ S(\ph'')\circ h_i
=\textsum_{i,(\ph)}\ph'\circ(a_i)_L\circ h_i\circ S(\ph'')
=\ph\triangleright\pi(x).
\vadjust{\vskip-4pt}
$$
Hence the required action of $E$ on itself is exactly the adjoint one.
To obtain the second isomorphism in (iii) one has to use the $H$-module map
$\nu:H_{\ad}\otimes A\to\End A$ such that
$h\otimes a\mapsto h\circ a_L$ instead of $\pi$. The map $\nu$ is bijective
since $\nu=\pi\circ\om$ where $\om:H_{\ad}\otimes A\to A\otimes H_{\ad}$
is a linear bijection defined by $h\otimes a\mapsto\sum_{(h)}h'a\otimes h''$
with inverse $a\otimes h\mapsto\sum_{(h)}h''\otimes S^{-1}(h')a$.
\endproof

\Remark.
Observe the appearance of $(E^*)\op$ as an intermediate link in (i). In view
of Proposition 4.1 the second isomorphism in (i) can be rewritten as
$(A\otimes A\op)^H\cong(E^*)\op$. If $\sum_ja'_j\otimes b'_j\in(A\otimes A)^H$
corresponds to $\xi\in E^*$ under this isomorphism then, by (iii) of
Proposition 4.1, $\<\xi,\psi\>=\sum_ja'_j(\psi b'_j)$ for all $\psi\in E$. This
together with (3.2) shows that we recover the embedding of $(E^*)\op$ in
$A\otimes A\op$ that gives rise to the Miyashita-Ulbrich action of $E^*$,
that is, $c\xi=\sum_ja'_jcb'_j$ for all $c\in A$. One has by symmetry
$(A\otimes A\op)^E\cong(H^*)\op$.
\endremark

\proclaim
Proposition 4.4.
Suppose that $V\subset H_{\ad}$ is an $H$-submodule which generates $H$ as an
algebra, and let $\{e_i\}$ and $\{e_i^*\}$ be dual bases for $V$ and $V^*$.
Then $E$ is generated as an algebra by the operators
$$
\vadjust{\vskip-10pt}
\Psi_\tau=\textsum_i\pi\bigl(\tau(e_i^*)\otimes e_i\bigr),\qquad\tau\in D(\*V).\eqno(4.2)
$$
\endproclaim

\Proof.
By Proposition 4.3(iii) $\pi$ maps the algebra $B''=(A\otimes B\op)^H$ where
$B=H_{\ad}$ isomorphically onto $E$. By Proposition 4.1(iv) $B''$ is generated by
$(A\otimes V)^H$. Finally, $(A\otimes V)^H\cong D(\*V)$ in such a way that
$\sum_i\tau(e_i^*)\otimes e_i$ corresponds to $\tau$.
\endproof

Let $\si\in D(V)$ and $\tau\in D(\*V)$ where $V\in H\mo$. Suppose that
$\{e_i\}$ and $\{e_i^*\}$ are dual bases for $V$ and $V^*$, respectively. Applying
$\tau\otimes\si$ to the $H$-invariant element $\sum_ie_i^*\otimes e_i\in
\*V\otimes V$, we get $\sum_i\tau(e_i^*)\otimes\si(e_i)\in(A\otimes A)^H$. Define
$\Up_{\tau\si}\in E^*$ by
$$
\vadjust{\vskip-6pt}
\<\Up_{\tau\si},\ph\>=\textsum_i\tau(e_i^*)\cdot\ph\si(e_i)
\qquad{\fam0for\ }\ph\in E\eqno(4.3)
$$
(note that the right hand side belongs to  $A^H=k$). Suppose that
$\{\tau_j\}$ and $\{\si_j\}$ are dual bases of $D(V)$
and $D(\*V)$ with respect to the pairing $D(\*V)\times D(V)\to k$ of
Proposition 1.4(iv).

\proclaim
Proposition 4.5.
For $\si\in D(V)$, $\,\tau\in D(\*V)$ and $\ph\in E$ one has
\mlines
$$
\<\Up_{\tau\si},S(\ph)\>=\textsum_i\ph\tau(e_i^*)\cdot\si(e_i),\eqno(4.4)
$$$$
\ph\si=\textsum_j\<\Up_{\tau_j\si},\ph\>\si_j,\qquad\qquad
\ph\tau=\textsum_j\<\Up_{\tau\si_j},S(\ph)\>\tau_j,\eqno(4.5)
$$$$
\ep(\Up_{\tau\si})=\<\tau,\si\>,\qquad\qquad
\De(\Up_{\tau\si})=\textsum_j\Up_{\tau\si_j}\otimes\Up_{\tau_j\si}.\eqno(4.6)
$$
\mlines
\endproclaim

\vskip-12pt
\Proof.
One rewrites the right hand side of (4.4) as
$$
\ph'\bigl(\textsum_{i,(\ph)}\tau(e_i^*)\cdot S(\ph'')\si(e_i)\bigr)
=\textsum_{i,(\ph)}\<\Up_{\tau\si},S(\ph'')\>\ph'(1),
\vadjust{\vskip-6pt}
$$
which collapses to $\<\Up_{\tau\si},S(\ph)\>$ since $\ph'(1)=\ep(\ph')$.
Using Proposition 1.4(v) with $\*V$ in place of $V$, we obtain
\mlines
$$
\ph\si(v)=\textsum_i\<v,e_i^*\>\ph\si(e_i)
=\textsum_{i,j}\si_j(v)\tau_j(e_i^*)\ph\si(e_i)
=\textsum_j\<\Up_{\tau_j\si},\ph\>\si_j(v),
$$$$
\ph\tau(\xi)
=\textsum_i\<e_i,\xi\>\ph\tau(e_i^*)
=\textsum_{i,j}\ph\tau(e_i^*)\si_j(e_i)\tau_j(\xi)
=\textsum_j\<\Up_{\tau\si_j},S(\ph)\>\tau_j(\xi)
\vadjust{\vskip-6pt}
$$
\mlines
for all $v\in V$ and $\xi\in V^*$, whence (4.5). Substituting $\ph=\id$ in
(4.3) yields the first formula in (4.6). The last formula follows from the
computation below, where we take $\ph,\psi\in E$ and use (4.3) and (4.5):
$$
\<\Up_{\tau\si},\ph\psi\>
=\textsum_i\tau(e_i^*)\ph\psi\si(e_i)
=\textsum_{i,j}\<\Up_{\tau_j\si},\psi\>\tau(e_i^*)\ph\si_j(e_i)
=\textsum_j\<\Up_{\tau\si_j},\ph\>\<\Up_{\tau_j\si},\psi\>.
$$

\proclaim
Proposition 4.6.
If $V\subset A$ is an $H$-invariant subspace which generates the algebra $A$ and
$\io:V\to A$ is the inclusion map, then $E^*$ is generated as an algebra by the
elements $\Up_{\tau\io}$ with $\tau\in D(\*V)$.
\endproclaim

\Proof.
By Proposition 4.3 the algebra $A''=(A\otimes A\op)^H$ is isomorphic to $E^*$. Under
this isomorphism the invariant $\sum_i\tau(e_i^*)\otimes\si(e_i)$ goes to
$\Up_{\tau\si}\circ S$ (see the remark following Proposition 4.3). Take $\si=\io$ and
apply Proposition 4.1(iv).
\endproof

\Remark.
One can generalize the statement above as follows. If $V\in H\mo$ and $\si\in D(V)$
is such that $\si(V)$ generates the algebra $A$, then the set
$\{\Up_{\tau\si}\mid\tau\in D(\*V)\}$ generates the algebra $E^*$. Quite similarly,
$\{\Up_{\tau\si}\mid\si\in D(V)\}$ generates $E^*$ as long as $\tau(\*V)$ generates $A$.
\endremark

Suppose next that $H$ is cocommutative. In this case $\*V=V^*$ for all $V\in H\mo$.
Given $\si\in D(V)$ and $\tau\in D(V^*)$, define $\Phi_{\tau\si}:A\to A$ by the rule
$$
\vadjust{\vskip-6pt}
\Phi_{\tau\si}(a)=\textsum_i\tau(e^*_i)a\si(e_i)\qquad{\fam0for\ }a\in A.\eqno(4.7)
$$
Let $\R=\sum_m\,\rho'_m\otimes\rho''_m$ be the triangular structure on $E$
corresponding to $R=1\otimes1$. Applying (2.4), we get
$\Phi_{\tau\si}(a)=\sum_{i,m}\tau(e^*_i)\rho''_m\si(e_i)\rho'_m(a)$.
Thus $\Phi_{\tau\si}=f(\Up_{\tau\si})$ where $f:E^*\to E$ is the map defined
by $\xi\mapsto\sum_m\<\xi,\rho''_m\>\rho'_m$. In particular,
$\Phi_{\tau\si}\in E$. Since $f$ is a homomorphism of coalgebras, (4.6)
translates into
$$
\ep(\Phi_{\tau\si})=\<\tau,\si\>,\qquad\qquad
\De(\Phi_{\tau\si})=\textsum_j\Phi_{\tau\si_j}\otimes\Phi_{\tau_j\si}.\eqno(4.8)
$$

\proclaim
Proposition 4.7.
Suppose that $H$ is cocommutative and $A$ is central simple. Denote by $C_V\subset E$
the subcoalgebra spanned by $\{\Phi_{\tau\si}\mid\si\in D(V),\tau\in D(V^*)\}$.
Then the assignment $V\mapsto C_V$ gives a bijection between the isomorphism
classes of irreducible $H$-modules and the simple subcoalgebras of $E$.
If $k$ is algebraically closed then $\dim C_V=(\dim V)^2$ for each irreducible
$V$.
\endproclaim

\Proof.
The irreducible left $H$-modules are in a bijective correspondence with the
irreducible left $E$-modules under the duality $D$. The latter correspond
bijectively to simple right $E^*$-comodules, and those to simple subcoalgebras
of $E^*$. By (4.5) the comodule structure on $D(V)$ is given by
$\si\mapsto\sum_j\si_j\otimes\Up_{\tau_j\si}$. It follows that the subcoalgebra
$C^V\subset E^*$ corresponding to $D(V)$ is spanned by $\Up_{\tau\si}$ with
$\si\in\nobreak D(V)$ and $\tau\in D(V^*)$. It is therefore simple when $V$ is irreducible.
Since $\rk\R=\dim E$ by Proposition 2.3, the map $f:E^*\to E$ is bijective. It is
therefore a coalgebra isomorphism, and so it maps simple subcoalgebras of $E^*$ to
simple subcoalgebras of $E$. Observe that $f(C^V)=C_V$. Finally,
$\dim C^V=(\dim D(V))^2=(\dim V)^2$ when $k$ is algebraically closed.
\endproof

\proclaim
Proposition 4.8.
Suppose that $H$ is cocommutative and $A$ is central simple.
If $V\subset A$ is an $H$-invariant subspace which generates the algebra $A$ and
$\io:V\to A$ is the inclusion map, then $E$ is generated as an algebra by the
operators $\Phi_{\tau\io}$ with $\tau\in D(V^*)$.
\endproclaim

\Proof.
Since $f:E^*\to E$ is bijective, it is an antiisomorphism of algebras.
So the statement follows at once from Proposition 4.6.
\endproof

\proclaim
Proposition 4.9.
Suppose that $H$ is cocommutative.
If $\c$ is the braiding of $E\mo$ determined by $\R$, then
\mlines
$$
\c_{D(V)D(V)}(\si\otimes\si')=\textsum_j\si_j\otimes\Phi_{\tau_j\si'}\si
\qquad{\fam0for\ }\si,\si'\in D(V),\eqno(4.9)
$$$$
\c_{D(V^*)D(V^*)}(\tau\otimes\tau')=\textsum_j\Phi_{\tau\si_j}\tau'\otimes\tau_j
\qquad{\fam0for\ }\tau,\tau'\in D(V^*).\eqno(4.10)
$$
\mlines
\endproclaim

\vskip-12pt
\Proof.
Recall that $D(V)\otimes D(V)\cong D(V\otimes V)$ as described in
Proposition 1.4(iii). So we have to evaluate both sides of (4.9) at $u\otimes v$
where $u,v\in V$. Now
$$
\si(v)\si'(u)=\textsum_i\<u,e_i^*\>\si(v)\si'(e_i)
=\textsum_{i,j}\si_j(u)\tau_j(e_i^*)\si(v)\si'(e_i)
\vadjust{\vskip-6pt}
$$
by Proposition 1.4(v). In view of (2.5) this verifies identity (4.9) since
$c_{VV}:V\otimes V\to V\otimes V$ is just the exchange of tensorands.
Similarly, (4.10) is obtained by taking $\xi,\eta\in V^*$ and checking that
$$
\tau(\eta)\tau'(\xi)=\textsum_i\<e_i,\eta\>\tau(e_i^*)\tau'(\xi)
=\textsum_{i,j}\tau(e_i^*)\tau'(\xi)\si_j(e_i)\tau_j(\eta).
\vadjust{\vskip-20pt}
$$
\hfill
\endproof

Using the bijection $f:E^*\to E$ we can transfer to $E$ the two $E$-module
structures on $E^*$ given by
$\ph\rightharpoonup\xi=\sum_{(\xi)}\<\xi'',\ph\>\xi'$ and
$\ph\rightharpoondown\xi=\sum_{(\xi)}\<\xi',S(\ph)\>\xi''$ where $\ph\in
E$, $\,\xi\in E^*$ and $\De\xi=\sum_{(\xi)}\xi'\otimes\xi''$. Denote the $E$-module structures on $E$ obtained in
this way likewise by $\rightharpoonup$ and $\rightharpoondown$. Thus
$$
\ph\rightharpoonup f(\xi)=f(\ph\rightharpoonup\xi)\qquad{\fam0and}\qquad
\ph\rightharpoondown f(\xi)=f(\ph\rightharpoondown\xi)\eqno(4.11)
$$

\proclaim
Proposition 4.10.
Suppose that $H$ is cocommutative and $A$ is central simple.
The action $\rightharpoondown$ makes $E$ into a left $E$-module algebra.
Given $\si\in D(V)$ and $\tau\in D(V^*)$ where $V\in H\mo$, one has
$$
\ph\rightharpoonup\Phi_{\tau\si}=\Phi_{\ph\tau\!,\,\si}\qquad
{\fam0and}\qquad
\ph\rightharpoondown\Phi_{\tau\si}=\Phi_{\tau\!,\,\ph\si}.\eqno(4.12)
$$
\endproclaim

\vskip-6pt
\Proof.
The first assertion follows from the fact that $(E^*)\op$ is  a left
$E$-module algebra with respect to $\rightharpoondown$ and $f$ is an algebra
isomorphism of $(E^*)\op$ onto $E$. In view of (4.6) and (4.5)
\mlines
$$
\ph\rightharpoonup\Up_{\tau\si}
=\textsum_j\<\Up_{\tau_j\si},\ph\>\Up_{\tau\si_j}=\Up_{\tau\!,\,\ph\si},
$$$$
\ph\rightharpoondown\Up_{\tau\si}
=\textsum_j\<\Up_{\tau\si_j},S(\ph)\>\Up_{\tau_j\si}=\Up_{\ph\tau\!,\,\si}.
\vadjust{\vskip-6pt}
$$
\mlines
Applying $f$ and using (4.11), we deduce (4.12).
\endproof

\section
5. Quantum Lie algebras and enveloping algebras

Suppose that $R=\sum_i\,r'_i\otimes r''_i$ is a triangular structure on
$H$. Then for every $V\in H\mo$ and $n>0$ the symmetric group $\S_n$
operates on $V^{\otimes n}$ via $H$-module automorphisms. Precisely, if
$t_p\in\S_n$ is the transposition of $p$ and $p+1$ for some $0<p<n$, then
the action of $t_p$ in $V^{\otimes n}$ is given by the operator
$\id\otimes c_{VV}\otimes\id$ where $c_{VV}:V\otimes V\to V\otimes V$ is the
braiding map in $H\mo$ applied to the $p$th and $(p+1)$th tensorands of
$V^{\otimes n}$. For $s\in\S_n$ denote by $s_V$ the corresponding
transformation of $V^{\otimes n}$. We will now specialize the definition of
Lie algebras in symmetric tensor categories proposed by Manin \cite{Man, Ch.~12}:

\proclaim
Definition 5.1.
A Lie algebra in $H\mo$ is a pair $(L,\la)$ where $L$ is an object and
$\la:L\otimes L\to L$ a morphism in $H\mo$ satisfying

\item(i)
Quantum anticommutativity: $\la$ is zero on the subspace of
$\S_2$-fixed elements in $L\otimes L$;

\item(ii)
Quantum Jacobi identity: $\la\circ(\id\otimes\la)\circ(\id+z_L+z_L^2)=0$ on
$L^{\otimes3}$ where $z\in\S_3$ is the cyclic permutation
$(123)$.

\endproclaim

The permutation $z=(12)(23)$ operates in $L^{\otimes3}$ as
$(c_{LL}\otimes\id)(\id\otimes c_{LL})=c_{L\otimes L,L}$. By the
naturality of the braiding one has
$(\id\otimes\la)\circ z_L=c_{LL}\circ(\la\otimes\id)$. Since
$\id+z_L+z_L^2=z_L\circ(\id+z_L+z_L^2)$ and $\la\circ c_{LL}=-\la$ by the
anticommutativity, the quantum Jacobi identity admits an equivalent formulation:
$$
\la\circ(\la\otimes\id)\circ(\id+z_L+z_L^2)=0.
$$

By a {\it quantum Lie algebra} we mean a Lie algebra in $H\mo$ for an
arbitrary triangular Hopf algebra, and we use brackets to denote the
multiplication in quantum Lie algebras. A relevant notion due to Majid
\cite{Maj94} does not require the tensor category to be symmetric but
involves more axioms. The morphisms of Lie algebras in $H\mo$ are defined in
an obvious way. Any $B\in H\alg$ can be made into a quantum Lie algebra by
means of the operation
$$
[b,d]_q=bd-\textsum_i\,(r''_id)(r'_ib)\qquad{\fam0for\ }b,d\in B.\eqno(5.1)
$$

\proclaim
Definition 5.2.
An enveloping algebra of a Lie algebra $L$ in $H\mo$ is a pair $(B,\io)$
where $B\in H\alg$ and $\io:L\to B$ is an embedding of Lie algebras in
$H\mo$ such that $\io(L)$ generates $B$ as an associative algebra.
\endproclaim

\proclaim
Lemma 5.3.
Suppose that $(B,\io)$ is an enveloping algebra of $L$. Let $\{e_i\}$ be a
basis for $L$ and $x_i=\io(e_i)$ for each $i$. Suppose that
$c_{LL}(e_i\otimes e_j)=\sum_{i'\!,j'}R_{ij}^{i'\!j'}e_{j'}\otimes e_{i'}$
and $[e_i,e_j]=\sum_l\ga_{ij}^le_l$ where $R_{ij}^{i'\!j'}\!,\,\ga_{ij}^l\in
k$. Then
$$
x_ix_j-\textsum_{i'\!,j'}R_{ij}^{i'\!j'}x_{j'}x_{i'}=\textsum_l\ga_{ij}^lx_l\,.\eqno(5.2)
$$
\endproclaim

\vskip-6pt
\Proof.
We have $c_{BB}(x_i\otimes x_j)=\sum_{i'\!,j'}R_{ij}^{i'\!j'}x_{j'}\otimes
x_{i'}$ since $\io$ is a morphism in $H\mo$. This shows that the left hand
side of (5.2) coincides with the quantum Lie product $[x_i,x_j]_q$ defined in
(5.1). Now (5.2) just expresses the property that $\io$ is a morphism
of quantum Lie algebras.
\endproof

Let us return to the settings (1.7). Suppose that $R$ and $\R$ are mutually
corresponding triangular structures on $H$ and $E=\End_HA$.

\proclaim
Lemma 5.4.
Let $V\in H\mo$ and $s\in\S_n$. Denote by $u\in\S_n$ the
involution that reverses the order of $1,\ldots,n$.

\item(i)
Under the isomorphism $D(V^{\otimes n})\cong D(V)^{\otimes n}$ of
Proposition {\rm1.4(iii)} $D(s_V)$ corresponds to $s_{D(V)}^{-1}$.

\item(ii)
Under the canonical isomorphisms $(V^{\otimes n})^*\cong(V^*)^{\otimes n}$
and $\*(V^{\otimes n})\cong(\*V)^{\otimes n}$ in $H\mo$ the operator $s_V^*$
corresponds to $(us^{-1}u^{-1})_{V^*}$ and $\*\!s_V$ to
$(us^{-1}u^{-1})_{\*\!V}$.

\endproclaim

\Proof.
It suffices to consider only transpositions $s=t_p$ (since $\S_n$ is generated
by those). In case (i) the assertion follows from the commutativity of diagram
(2.5). For left duals we use the diagram
$$
\diagram{
V^*\otimes W^*&{}\cong{}&(W\otimes V)^*\cr
\noalign{\smallskip}
\mapdown{c_{V^*W^*}}{}&&\mapdown{}{c_{VW}^*}\cr
\noalign{\smallskip}
W^*\otimes V^*&{}\cong{}&(V\otimes W)^*\cr
}
$$
where $V,W\in H\mo$. The commutativity of this diagram is a general property
of braided tensor categories. In case of Hopf algebras this property reduces
to the identity
$\sum_i\,\<\xi,S(r'_i)v\>\<\eta,S(r''_i)w\>=\sum_i\,\<\xi,r'_iv\>\<\eta,r''_iw\>$
where $\xi\in V^*$, $\eta\in W^*$, $v\in V$, $w\in W$. The latter does hold
since $(S\otimes S)(R)=R$ \cite{Dr89a}. The isomorphisms in (ii) reverse the
tensorands, and therefore $(t_p)_V^*$ corresponds to $(t_{n-p-1})_{V^*}$.
The case of right duals is similar.
\endproof

\proclaim
Proposition 5.5.
Suppose that $(L,\la)$ is a Lie algebra in $H\mo$. Then:

\item(i)
$\bigl(D(L^*),D(\la^*)\bigr)$ and $\bigl(D(\*\!L),D(\*\!\la)\bigr)$ are Lie algebras
in $E\mo$.

\item(ii)
If $(B,\io)$ is an enveloping algebra of $L$, then
$\bigl(D\bigl((B\op)^*\bigr),D(\io^*)\bigr)$ is an enveloping algebra of
$D(L^*)$ and $\bigl(D\bigl(\*(B\op)\bigr),D(\*\!\io)\bigr)$ that of
$D(\*\!L)$.

\item(iii)
If $A$ {\rm(}respectively, $H^*_\rightharpoonup$, $H_{\ad}${\rm)} is an enveloping
algebra of $L$, then $E^*_\rightharpoonup$ {\rm(}respectively, $A$,  $E_{\ad}${\rm)}
is an enveloping algebra of both $D(L^*)$ and $D(\*\!L)$.

\endproclaim

\Proof.
Let $\la'=D(\la^*)$ and $L'=D(L^*)$. The fixed point subspace $I=(L\otimes L)^{\S_2}$ can be
characterized as the kernel of the transformation $\id-t_L$ where $t\in\S_2$ is
the transposition. Since the functor $D(?^*)$ is exact, it follows from Lemma 5.4
that $D(I^*)$ coincides with the subspace of $\S_2$-fixed elements in
$L'\otimes L'$. Clearly $\la'$ is zero on $D(I^*)$. Since $\{e,z,z^2\}$ is a
normal subgroup of $\S_3$, the functor $D(?^*)$ transforms the quantum Jacobi
identity for $L$ into the identity $\la'\circ(\la'\otimes\id)\circ(\id+z_{L'}+z_{L'}^2)=0$.

Let $\mu:B\otimes B\to B$ denote the associative multiplication. Put $B'=D(B^*)$ and
$\mu'=D(\mu^*)$. The quantum Lie multiplication in $B$ is $\mu-\mu\circ c_{BB}$ and that
in $B'$ is $\mu'-\mu'\circ\c_{B'B'}$. These two correspond to one another
under the functor $D(?^*)$. By functoriality $\io'=D(\io^*)$ is a homomorphism of
Lie algebras in $E\mo$. It is an embedding since $D$ is exact. The image of
$\io'$
generates $B'$ by Proposition 4.1(iv).
Right duals are treated similarly. Part (iii) follows from Proposition 4.3.
\endproof

\section
6. Reduced enveloping algebras as Galois algebras

Let $k$ be a field of characteristic $p>0$ and $\g$ a finite dimensional
$p$-Lie algebra over $k$. For $\xi\in\g^*$ denote by $U_\xi(\g)$ the
corresponding reduced enveloping algebra of $\g$. Thus $U_\xi(\g)$ is the
factor algebra of the universal enveloping algebra $U(\g)$ by its ideal generated by
central elements $x^p-x^{[p]}-\xi(x)^{p}1$ with $x\in\g$. We identify $\g$
with a Lie subalgebra in each $U_\xi(\g)$. There is an algebra homomorphism
$$
\De:U_\xi(\g)\to U_0(\g)\otimes U_\xi(\g),\qquad
x\mapsto1\otimes x+x\otimes1{\fam0\ for\ }x\in\g.\eqno(6.1)
$$
In particular, $\De$ makes $U_0(\g)$ into a cocommutative Hopf algebra, and
each $U_\xi(\g)$ is a $U_0(\g)$-comodule algebra. Next, the adjoint action
of $\g$ on itself extends to a restricted action of $\g$ by derivations
on $U_\xi(\g)$. In general, if $\g$ operates as a $p$-Lie algebra of
derivations on an algebra $A$, one may regard $A$ as a left $U_0(\g)$-module
algebra. We say that $A$ is {\it $\g$-Galois\/} if $A$ is $U_0(\g)^*$-Galois
with respect to the corresponding comodule structure.

\proclaim
Proposition 6.1.
The algebra $U_\xi(\g)$ is $\g$-Galois if and only if $U_\xi(\g)$ is
central simple. In this case the comodule structure {\rm(6.1)} corresponds to the
Miyashita-Ulbrich action of $U_0(\g)^*$.
\endproclaim

\Proof.
Let $H=U_0(\g)$ and $A=U_\xi(\g)$. There is then an isomorphism of algebras
$\nu:A\#H\to A\otimes A\op$ defined by $\nu(u\#1)=u\otimes1$ and
$\nu(1\#x)=x\otimes1-1\otimes x$ for $u\in A$ and $x\in\g$ (one can easily
construct the inverse). If $\ka:A\otimes A\op\to\End A$ is the canonical
homomorphism given by left and right multiplications, then $\ka\circ\nu$
coincides with (1.3). Hence $\ka$ is bijective if and only if so is map (1.3).
This proves the first assertion. As follows from (3.1), the
comodule structure $b\mapsto\sum_i h_i\otimes bh_i^*$ corresponding to the
Miyashita-Ulbrich action is characterized by the condition that
$\pi'(\sum_i bh_i^*\otimes h_i)$, where $\pi'$ is defined in (1.4), is the
left multiplication by $b$ in $A$.
Since $xu=(\ad x)u+ux$, we see that the left multiplication by $x\in\g$ coincides
with $\pi'(1\otimes x+x\otimes1)$. Hence (6.1) gives this comodule structure
on elements of $\g$. Since $U_\xi(\g)$ is generated by $\g$, formula (6.1)
determines the structure of a comodule algebra uniquely.
\endproof

\Remark.
The algebra $U_\xi(\g)$ is $U_0(\g)$-Galois with respect to (6.1) for any
$\g$ and $\xi$. Furthermore, the corresponding Miyashita-Ulbrich action is
exactly the adjoint action of $U_0(\g)$ on $U_\xi(\g)$. Hence Proposition 6.1 is
an instance of general results in \cite{Ul82} and \cite{Doi89} which were
reproduced in Proposition 3.3.
\endremark

If $\g$ is Frobenius then $\dim\g$ is even and the maximum dimension of irreducible
$\g$-modules over an algebraically closed extension field of $k$ equals
$p^{{1\over2}\dim\g}$ by \cite{Mil75} or \cite{Pre99, Th.~4.4}. Note that
$\dim U_\eta(\g)=p^{\dim\g}$ for every $\eta\in\g^*$ by the PBW theorem. In case of an algebraically
closed $k$ this means that the algebras $U_\eta(\g)$ are simple for
all $\eta$ in a nonempty Zariski open subset of $\g^*$. This provides a
plentiful supply of central simple $\g$-Galois algebras. Proposition 6.2 below adapts to
characteristic $p$ one of the arguments due to Ooms \cite{Oo80, Th.~3.3} who studied Frobenius
Lie algebras in characteristic 0. Recall that $\g$ is {\it unimodular} if $\ad_\g x$ has trace $0$
for every $x\in\g$.

\proclaim
Proposition 6.2.
If $\g$ is Frobenius and unimodular, then\/ $\dim\g\equiv0\mod{2p}$.
\endproclaim

\Proof.
Let $\xi\in\g^*$ be a linear function such that the associated alternating bilinear form $\be_\xi$
is nondegenerate. Then the map $\g\to\g^*$ given by $x\mapsto\xi\circ(\ad x)$ is bijective. There
exists $x\in\g$ such that $\xi\circ(\ad x)=\xi$, that is, $\xi([xy])=\xi(y)$ for all
$y\in\g$. Consider the weight space decomposition
$\g=\oplus\g_\la$ with respect to the adjoint transformation $\ad x$. Here
$$
\g_\la=\{y\in\g\mid(\ad x-\la\id)^ny=0\}
$$
for $\la\in k$ and sufficiently big $n$.
Our choice of $x$ ensures that $\xi(\g_\la)=0$ whenever $\la\ne1$, and so
$\xi([\g_\la\g_\mu])=0$ for all $\la,\mu\in k$ such that $\la+\mu\ne1$. The
restriction of $\be_\xi$ induces a nondegenerate pairing $\g_\la\times\g_{1-\la}\to k$
for each $\la\in k$. This shows that $\dim\g_\la=\dim\g_{1-\la}$. If $2\la=1$
then $\dim\g_\la$ is even as $\g_\la$ admits a nondegenerate alternating
bilinear form. Hence $\tr(\ad x)=\sum_\la(\dim\g_\la)\la={\dim\g\over2}1$.
If $\tr(\ad x)=0$, then ${\dim\g\over2}\equiv0\mod p$.
\endproof

Further on we assume that
$U_\xi(\g)$ is central simple and $E_\xi=\End_\g U_\xi(\g)$.

\proclaim
Theorem 6.3.
Under the above assumption one has:

\item(i)
There is a unique comultiplication on $E_\xi$ with respect to which $E_\xi$
is a Hopf algebra and $U_\xi(\g)$ is a left $E_\xi$-module algebra.

\item(ii)
There is a triangular structure $\R=\sum_m\rho'_m\otimes\rho''_m$ of maximal
rank  on $E_\xi$ such that $uv=\sum_m(\rho''_mv)(\rho'_mu)$ for all
$u,v\in U_\xi(\g)$. In particular, $E_\xi^*\cong E_\xi\op$ as Hopf algebras.

\item(iii)
The antipode $S$ of $E_\xi$ has order $2p$ or $2$ {\rm(}with the exception of the
case $\g=0$ when $S=\id${\rm)}. Moreover, $S^2=\id$ if and only if $\g$ is unimodular.

\endproclaim

\Proof.
Assertions (i) and (ii) are special cases of Theorem 1.6, Propositions 2.2 and
2.3. The order of the antipode is determined in Theorem 3.8. The Hopf
algebra $U_0(\g)$ is generated by primitive elements, and so it is pointed irreducible.
In particular, it contains a single grouplike element. The modular function
$\al:U_0(\g)\to k$ is determined by $\al(x)=\tr(\ad_\g x)$ for $x\in\g$
(cf.\ \cite{Jan, Remark in I.9.7}). In particular, $U_0(\g)$ is
unimodular if and only if so is $\g$ \cite{Lar69, p.~91}. According to (6.1)
the Miyashita-Ulbrich action of $\al$ on $U_\xi(\g)$ is given by
$x\al=x+\al(x)1$ for $x\in\g$. Hence $x\al^p=x+p\al(x)1=x$ for all $x$.
Since $\al$ operates as an automorphism, we have $\al^p=1$. Thus the order
of $\al$ is either $p$ or $1$. Finally, note that the order of $S$ is even
unless $E_\xi$, hence also $U_0(\g)$ are commutative, which is only possible
when $\g=0$.
\endproof

As explained in Proposition 1.4, there is a duality $D$ between $U_0(\g)\mo$ and
$E_\xi\mo$. We will use the $E_\xi$-module $D(V)=\Hom_\g\bigl(V,U_\xi(\g)\bigr)$
when $V$ is either the adjoint $\g$-module $\g$ or the coadjoint one
$\g^*$. Let $\{e_i\}$ and $\{e_i^*\}$ be dual bases for $\g$ and $\g^*$.
Recall that the nondegenerate pairing $D(\g^*)\times D(\g)\to k$ is defined by the
formula $\<\tau,\si\>=\sum_i\tau(e_i^*)\si(e_i)$ for $\tau\in D(\g^*)$ and
$\si\in D(\g)$. Both Proposition 4.4 and Proposition 4.8 describe generating subspaces for
the algebra $E_\xi$. We may take $V=\g$ and $\io:\g\to U_\xi(\g)$ the canonical
embedding. Then $\Psi_\tau,\Phi_{\tau\io}\in E_\xi$ are defined for all
$\tau\in D(\g^*)$. As
$$
\Psi_\tau u=\textsum_i\tau(e_i^*)(\ad e_i)u=\textsum_i\tau(e_i^*)e_iu
-\textsum_i\tau(e_i^*)ue_i=\<\tau,\io\>u-\Phi_{\tau\io}u
\vadjust{\vskip-6pt}
$$
for all $u\in U_\xi(\g)$, one sees that $\Psi_\tau=\<\tau,\io\>\id-\Phi_{\tau\io}$.
Consider the triangular structure $\R$ on $E_\xi$ corresponding to
the triangular structure $R=1\otimes1$ on $U_0(\g)$.

\proclaim
Theorem 6.4.
$D(\g^*)$ is a Lie algebra in $E_\xi\mo$ with
respect to the multiplication
$$
\vadjust{\vskip-6pt}
[\tau,\tau'](\xi)=\textsum_{i,j}\<\xi,[e_j,e_i]\>\tau(e_i^*)\tau'(e_j^*),
\qquad\tau,\tau'\in D(\g^*),\ \xi\in\g^*.\eqno(6.2)
$$
The algebra $E_\xi$ equipped with either the adjoint action on itself or the
action $\rightharpoondown$ described in Proposition {\rm4.10} is an enveloping algebra
of $D(\g^*)$. The corresponding embeddings $D(\g^*)\to E_\xi$ are given by
$\tau\mapsto\Psi_\tau$ and $\tau\mapsto\Phi_{\tau\io}$, respectively.
\endproclaim

\Proof.
Equipped with the adjoint action on itself, $\g$ is a Lie algebra in
$U_0(\g)\mo$. Indeed, the quantum anticommutativity and Jacobi identity in
$U_0(\g)\mo$ coincide with the their ordinary versions since $R=1\otimes1$.
Now $D(\g^*)$ is a quantum Lie algebra by Proposition 5.5. Next, $U_\eta(\g)$ is an
enveloping algebra of $\g$ for any $\eta\in\g^*$. By Proposition 5.5 the functor
$D(?^*)$ transforms $U_\eta(\g)$ into an enveloping algebra of $D(\g^*)$. If
$\eta=0$ (respectively, $\eta=\xi$), the result is $(E_\xi)_{\ad}$ (respectively,
$(E_\xi^*)_\rightharpoonup$). The $E_\xi$-module algebra in the second case is
isomorphic to $(E_\xi^*)\op$ with the action of $E_\xi$ given by
$\rightharpoondown$. Finally, $(E_\xi^*)\op\cong E_\xi$ as $E_\xi$-module
algebras with respect to $\rightharpoondown$. The embedding of $D(\g^*)$ into
$(E_\xi)_{\ad}$ is described in Proposition 4.4 with $V=\g$. The embedding of
$D(\g^*)$ into $(E_\xi^*)\op$ is given by $\tau\mapsto\Up_{\tau\io}$
(cf.\ Proposition 4.6). The element $\Up_{\tau\io}$ corresponds to $\Phi_{\tau\io}$
in $E_\xi$.
\endproof

\Example 1.
Let $\g$ be the 2-dimensional Lie algebra with a basis $e_0,e_1$ and the
multiplication $[e_0,e_1]=e_1$. The $[p]$-map on $\g$ is given by
$e_0^{[p]}=e_0$ and $e_1^{[p]}=0$. If $\xi(e_1)\ne0$ for $\xi\in\g^*$, then
$U_\xi(\g)$ has a $p$-dimensional absolutely irreducible module \cite{St,
6.9}. Since $\dim U_\xi(\g)=p^2$, this means that $U_\xi(\g)$ is central
simple. One can easily check that the linear functions $\xi$ satisfying
$\xi(e_1)\ne0$ constitute a conjugacy class with respect to the automorphism
group of $\g$. We may assume therefore that $\xi(e_0)=0$ and $\xi(e_1)=1$.
Then $e_0^p=e_0$ and $e_1^p=1$ in $U_\xi(\g)$.

The tables below give dual bases for $D(\g^*)$ and $D(\g)$ with respect to
the pairing $D(\g^*)\times D(\g)\to k$, the multiplication in the quantum Lie
algebra $D(\g^*)$ and the action of elements $\ph_i=\Phi_{\tau_i\si_1}$
in the $E_\xi$-module $D(\g^*)$:
$$
\table3
&e_0^*&e_1^*\cr
\tau_0&1&-e_0e_1^{-1}\cr
\tau_1&0&e_1^{-1}\cr
\endtable\qquad
\table3
&e_0&e_1\cr
\si_0&1&0\cr
\si_1&e_0&e_1\cr
\endtable\qquad
\table3
&\tau_0&\tau_1\cr
\tau_0&0&-\tau_1\cr
\tau_1&\tau_1&0\cr
\endtable\qquad
\table3
&\tau_0&\tau_1\cr
\ph_0&0&\tau_1\cr
\ph_1&\tau_0-\tau_1&\tau_1\cr
\endtable
$$
One sees that $D(\g^*)$ is in fact an ordinary Lie algebra isomorphic to
$\g$. The algebra $E_\xi$ is generated by $\ph_0,\ph_1$ (note that
$\si_1$ is the canonical embedding of $\g$ into $U_\xi(\g)$). One checks
that $\Phi_{\tau_0\si_0}=\id$ and $\Phi_{\tau_1\si_0}=0$.
The braiding map $\c_{D(\g^*)D(\g^*)}$ is computed from (4.10) as
$$
\vbox{\halign{\hfil$#$&${}\mapsto#$\hfil&\qquad\qquad\hfil$#$&${}\mapsto#$\hfil\cr
\tau_0\otimes\tau_0&\tau_0\otimes\tau_0,
&\tau_0\otimes\tau_1&\tau_1\otimes(\tau_0+\tau_1),\cr
\noalign{\smallskip}
\tau_1\otimes\tau_0&(\tau_0-\tau_1)\otimes\tau_1,
&\tau_1\otimes\tau_1&\tau_1\otimes\tau_1.\cr
}}
$$
Only one of the quadratic relations (5.2) is independent in the present
settings. For $i=0$, $j=1$ it gives $\ph_0\ph_1-\ph_1(\ph_0+\ph_1)=-\ph_1$.
The $p$th powers of $\ph_0,\ph_1$ can be figured out by computing
straightforwardly in the algebra $U_\xi(\g)\otimes U_\xi(\g)\op$. Note that
$\ph_0,\ph_1$ are represented by the elements $1\otimes e_0-e_0e_1^{-1}\otimes e_1$
and $e_1^{-1}\otimes e_1$ in the latter algebra. The coproduct in $E_\xi$ is
derived from (4.8). Altogether we get the following relations:
\mlines
$$
[\ph_0,\ph_1]=\ph_1^2-\ph_1,
\qquad\quad\ph_0^p=\ph_0,\qquad\quad\ph_1^p=1,
$$$$
\De(\ph_0)=1\otimes\ph_0+\ph_0\otimes\ph_1,\qquad\quad
\De(\ph_1)=\ph_1\otimes\ph_1,
$$$$
\ep(\ph_0)=0,\qquad\quad\ep(\ph_1)=1,\qquad\quad
S(\ph_0)=-\ph_0\ph_1^{-1},\qquad\quad S(\ph_1)=\ph_1^{-1}.
$$
\mlines
In this example one can recognize the Hopf algebra $H_p$ described in
\cite{Rad77, p.~158}. This Hopf algebra appears also in the family of
Hopf algebras $L_n$, $n>0$, constructed in \cite{Taft80, Lemma 7}.
In general $\dim L_n=p^{2n}$ and the antipode of $L_n$ has order $2p^n$.
One can suspect that $L_n$ can be realized as the endomorphism Hopf algebra
of a Galois algebra for the $n$th Frobenius kernel of the 2-dimensional
nonabelian algebraic group.

\Example 2.
Let $\g$ be the 4-dimensional Lie algebra with a basis $e_0,e_1,e_2,e_3$ and
nonzero commutators
$$
[e_0,e_1]=ae_1,\qquad
[e_0,e_2]=be_2,\qquad
[e_0,e_3]=(a+b)e_3,\qquad
[e_1,e_2]=e_3
$$
where $a,b$ are integers such that $a+b\not\equiv0\mod p$.
One has $e_0^{[p]}=e_0$ and $e_i^{[p]}=0$ for $i>0$. This Lie algebra is
unimodular only when $p=2$. The representations of $\g$ can be easily
determined using methods of \cite{St78}, \cite{Wei71}. If $\xi(e_3)\ne0$
then $U_\xi(\g)$ has an absolutely irreducible module of dimension $p^2$
which is induced from a onedimensional module over the subalgebra of $\g$
spanned by $e_2,e_3$.  Hence $U_\xi(\g)$ is central simple in this case. The
automorphism group operates transitively on the set of linear functions
satisfying $\xi(e_3)\ne0$. So we may assume that $\xi(e_3)=1$ and
$\xi(e_i)=0$ for $i\ne3$. The data comprised in the tables below has the
same meaning as in the previous example:
\mlines
$$
\table5
&e_0^*&e_1^*&e_2^*&e_3^*\cr
\tau_0&1&ae_2e_3^{-1}&-be_1e_3^{-1}&(be_1e_2-ae_2e_1-e_0e_3)e_3^{-2}\cr
\tau_1&0&e_3^d&0&-e_1e_3^{d-1}\cr
\tau_2&0&0&e_3^c&-e_2e_3^{c-1}\cr
\tau_3&0&0&0&e_3^{-1}\cr
\endtable
\vadjust{\medskip}
$$$$
\table5
&e_0&e_1&e_2&e_3\cr
\si_0&1&0&0&0\cr
\si_1&-ae_2e_3^c&e_3^{c+1}&0&0\cr
\si_2&be_1e_3^d&0&e_3^{d+1}&0\cr
\si_3&e_0&e_1&e_2&e_3\cr
\endtable\qquad
\table5
&\tau_0&\tau_1&\tau_2&\tau_3\cr
\tau_0&-ab\tau_3&-a\tau_1&-b\tau_2&-(a+b)\tau_3\cr
\tau_1&a\tau_1&0&-\tau_3&0\cr
\tau_2&b\tau_2&\tau_3&0&0\cr
\tau_3&(a+b)\tau_3&0&0&0\cr
\endtable
$$
\mlines
Here $c,d$ are integers such that $b+c(a+b)\equiv0$, $\,a+d(a+b)\equiv0$,
hence also $c+d+1\equiv0\mod p$. The multiplication in $D(\g^*)$ is no longer
anticommutative in the ordinary sense as $[\tau_0\tau_0]\ne0$. The elements
$\ph_i=\Phi_{\tau_i\si_3}$ generate $E_\xi$, and the relations fulfilled in
$E_\xi$ are listed below. We omit the details of the computations:
\mlines
$$
[\ph_0,\ph_1]=a\ph_1(2\ph_3-1),\qquad
[\ph_0,\ph_2]=b\ph_2(2\ph_3-1),\qquad
[\ph_0,\ph_3]=(a+b)(\ph_3^2-\ph_3),
$$$$
[\ph_1,\ph_2]=\ph_3^2-\ph_3,\qquad
[\ph_1,\ph_3]=[\ph_2,\ph_3]=0,
$$$$
\ph_0^p=\ph_0,\qquad\ph_1^p=\ph_2^p=0,\qquad\ph_3^p=1,
$$$$
\De(\ph_0)=1\otimes\ph_0-a\ph_2\ph_3^c\otimes\ph_1
+b\ph_1\ph_3^d\otimes\ph_2+\ph_0\otimes\ph_3,
$$$$
\De(\ph_1)=\ph_3^{-d}\otimes\ph_1+\ph_1\otimes\ph_3,\qquad
\De(\ph_2)=\ph_3^{-c}\otimes\ph_2+\ph_2\otimes\ph_3,\qquad
\De(\ph_3)=\ph_3\otimes\ph_3,
$$$$
\ep(\ph_0)=\ep(\ph_1)=\ep(\ph_2)=0,\qquad\ep(\ph_3)=1,
$$$$
S(\ph_0)=(b\ph_1\ph_2-a\ph_2\ph_1-\ph_0\ph_3)\ph_3^{-2},
$$$$
S(\ph_1)=-\ph_1\ph_3^{d-1},\qquad
S(\ph_2)=-\ph_2\ph_3^{c-1},\qquad
S(\ph_3)=\ph_3^{-1}.
$$
\mlines
The assignments $\ph_0\mapsto Z$, $\ph_1\mapsto X$, $\ph_2\mapsto Y$,
$\ph_3\mapsto1$ define a homomorphism of $E_\xi$ onto the Hopf algebra of
dimension $p^3$ with generators $X,Y,Z$ and relations
\mlines
$$
[Z,X]=aX,\qquad[Z,Y]=bY,\qquad[X,Y]=0,\qquad Z^p=Z,\qquad X^p=Y^p=0,
$$$$
\De(Z)=1\otimes Z-aY\otimes X+bX\otimes Y+Z\otimes1,
$$$$
\De(X)=1\otimes X+X\otimes1,\qquad\De(Y)=1\otimes Y+Y\otimes1.
$$
\mlines
It may be amusing to compare this Hopf algebra with the Hopf algebra $H$ from
\cite{Taft75}.

The Hopf algebras in examples 1, 2 are pointed. In the light of attempts to
classify pointed Hopf algebras (e.g., \cite{And01}) it may be of interest to single out those Hopf
algebras among the $E_\xi$'s.

\proclaim
Proposition 6.5.
Suppose that $k$ is algebraically closed and $U_\xi(\g)$ is simple. Then $E_\xi$ is pointed if and
only if $[\g,\g]$ consists of $[p]$-nilpotent elements.
\endproclaim

\Proof.
In order that $E_\xi$ be pointed it is necessary and sufficient that
$E_\xi^*$ have only onedimensional irreducible modules. Since
$E_\xi^*\cong E_\xi\cong U_0(\g)$ as algebras, an equivalent condition is
that $U_0(\g)$ has only onedimensional irreducible modules or, in other words,
$[\g,\g]$ annihilates the irreducible $U_0(\g)$-modules. In view of Engel's
theorem \cite{St, Ch.~1, Cor.~3.2} this can be rephrased by saying that all
elements of $[\g,\g]$ operate nilpotently in $U_0(\g)$-modules. Taking a
faithful $U_0(\g)$-module, we deduce that $x^{[p]^n}=0$ for all
$x\in[\g,\g]$ when $n$ is sufficiently big.
\endproof

\Remark.
If $k$ is not algebraically closed, one must add the condition that $\g$
contains a split torus of maximal dimension.
\endremark

\proclaim
Conjecture 6.6.
Let $\g$ be any $p$-Lie algebra and $\xi,\eta\in\g^*$.

\item(i)
$U_\xi(\g)$ is central simple if and only if the alternating bilinear
form $\be_\xi$ associated with $\xi$ is nondegenerate.

\item(ii)
The Hopf algebras $E_\xi$ and $E_\eta$ are isomorphic whenever both $U_\xi(\g)$
and $U_\eta(\g)$ are central simple.

\endproclaim

Currently I can prove that (i) holds if either $\g$ is solvable or
$k$ is algebraically closed and $\ad\g\subset\Lie G$ where $G$ is the
automorphism group of $\g$. In the second case (ii) is also fulfilled. These
results will appear elsewhere.

\references
\nextref
Ali01
\auth
E.,Aljadeff;P.,Etingof;S.,Gelaki;D.,Nikshych;
\endauth
\paper{On twisting of finite-dimensional Hopf algebras}
arXiv:math.QA/0107167.

\nextref
And01
\auth
N.,Andruskiewitsch;H.-J.,Schneider;
\endauth
\paper{Pointed Hopf algebras}\hfil\break\hbox{}\hfill
arXiv:math.QA/0110136.

\nextref
Bea
\auth
M.,Beattie;K.-H.,Ulbrich;
\endauth
\paper{A Skolem-Noether theorem for Hopf algebra actions}
\journal{Comm. Algebra}
\Vol{18}
\Year{1990}
\Pages{3713--3724}

\nextref
Bl89
\auth
R.J.,Blattner;S.,Montgomery;
\endauth
\paper{Crossed products and Galois extensions of Hopf algebras}
\journal{Pacific J. Math.}
\Vol{137}
\Year{1989}
\Pages{37--54}

\nextref
Chase69
\auth
S.U.,Chase;M.,Sweedler;
\endauth
\book{Hopf Algebras and Galois Theory}
\bookseries{Lecture Notes in Mathematics}
\Vol{97}
\publisher{Springer}
\Year{1969}

\nextref
Coh90
\auth
M.,Cohen;D.,Fischman;S.,Montgomery;
\endauth
\paper{Hopf Galois extensions, smash products, and Morita equivalence}
\journal{J. Algebra}
\Vol{133}
\Year{1990}
\Pages{351--372}

\nextref
Doi86
\auth
Y.,Doi;M.,Takeuchi;
\endauth
\paper{Cleft comodule algebras for a bialgebra}
\journal{Comm. Algebra}
\Vol{14}
\Year{1986}
\Pages{801--817}

\nextref
Doi89
\auth
Y.,Doi;M.,Takeuchi;
\endauth
\paper{Hopf-Galois extensions of algebras, the Miyashita-Ulbrich action, and Azumaya algebras}
\journal{J. Algebra}
\Vol{121}
\Year{1989}
\Pages{488--516}

\nextref
Dr87
\auth
V.G.,Drinfel'd;
\endauth
\paper{Quantum groups}
\journal{J. Sov. Math.}
\Vol{41}
\Year{1988}
\Pages{898-915}

\nextref
Dr89a
\auth
V.G.,Drinfel'd;
\endauth
\paper{On almost cocommutative Hopf algebras}
\journal{Leningrad Math. J.}
\Vol{1}
\Year{1990}
\Pages{321--342}

\nextref
Dr89b
\auth
V.G.,Drinfel'd;
\endauth
\paper{Quasi-Hopf algebras}
\journal{Leningrad Math. J.}
\Vol{1}
\Year{1990}
\Pages{1419-1457}

\nextref
Et00
\auth
P.,Etingof;S.,Gelaki;
\endauth
\paper{The classification of triangular semisimple and cosemisimple Hopf algebras over an algebraically closed field}
\journal{Internat. Math. Res. Notices}
\nombre{5}
\Year{2000}
\Pages{223--234}

\nextref
Fai
\auth
C.,Faith;
\endauth
\book{Algebra II Ring Theory}
\bookseries{Grundlehren der mathematischen Wissenschaften}
\Vol{191}
\publisher{Springer}
\Year{1976}

\nextref
Gr96
\auth
C.,Greither;
\endauth
\paper{On the transformation of the Hopf algebra in a Hopf Galois extension}
\journal{Comm. Algebra}
\Vol{24}
\Year{1996}
\Pages{737--747}

\nextref
Jac
\auth
N.,Jacobson;
\endauth
\book{Basic Algebra II}
\publisher{Freeman}
\Year{1980}

\nextref
Jan
\auth
J.,Jantzen;
\endauth
\book{Representations of algeraic groups}
\publisher{Academic Press}
\Year{1987}

\nextref
Joy93
\auth
A.,Joyal;R.,Street;
\endauth
\paper{Braided tensor categories}
\journal{Adv. Math.}
\Vol{102}
\Year{1993}
\Pages{20--78}

\nextref
Kas
\auth
C.,Kassel;
\endauth
\book{Quantum Groups}
\bookseries{Graduate Texts in Mathematics}
\Vol{155}
\publisher{Springer}
\Year{1995}

\nextref
Kop
\auth
M.,Koppinen;
\endauth
\paper{A Skolem-Noether theorem for coalgebra measurings}
\journal{Arch. Math.}
\Vol{57}
\Year{1991}
\Pages{34--40}

\nextref
Kr76
\auth
H.F.,Kreimer;P.M.I.,Cook;
\endauth
\paper{Galois theories and normal bases}
\journal{J. Algebra}
\Vol{43}
\Year{1976}
\Pages{115--121}

\nextref
Kr81
\auth
H.F.,Kreimer;M.,Takeuchi;
\endauth
\paper{Hopf algebras and Galois extensions of an algebra}
\journal{Indiana Univ. Math.~J.}
\Vol{30}
\Year{1981}
\Pages{675--692}

\nextref
Lar69
\auth
R.G.,Larson;M.E.,Sweedler;
\endauth
\paper{An associative orthogonal bilinear form for Hopf algebras}
\journal{Amer. J. Math.}
\Vol{91}
\Year{1969}
\Pages{75--94}

\nextref
Maj94
\auth
Sh.,Majid;
\endauth
\paper{Quantum and braided-Lie algebras}
\journal{J. Geom. Phys.}
\Vol{13}
\Year{1994}
\Pages{307-356}

\nextref
Maj
\auth
Sh.,Majid;
\endauth
\book{Foundations of Quantum Group Theory}
\publisher{Cambridge University Press}
\Year{1995}

\tracingmacros=1
\nextref
Man
\auth
Yu.,Manin;
\endauth
\book{Quantum Groups and Non-Commutative Geometry}
\bookseries{Publ. CRM}
\publisher{Universit\'e de Montr\'eal}
\Year{1988}

\nextref
Mas
\auth
A.,Masuoka;
\endauth
\paper{Coalgebra actions on Azumaya algebras}
\journal{Tsukuba J. Math.}
\Vol{14}
\Year{1990}
\Pages{107--112}

\nextref
Mil75
\auth
A.A.,Mil'ner;
\endauth
\paper{Irreducible representations of modular Lie algebras}
\journal{Math. USSR Izv.}
\Vol{9}
\Year{1975}
\Pages{1169--1187}

\nextref
Mo
\auth
S.,Montgomery;
\endauth
\book{Hopf algebras and Their Actions on Rings}
\bookseries{CBMS Regional Conference Series in Mathematics}
\Vol{82}
\publisher{American Mathematical Society}
\Year{1993}

\nextref
Nich89
\auth
W.D.,Nichols;M.B.,Zoeller;
\endauth
\paper{A Hopf algebra freeness theorem}
\journal{Amer. J. Math.}
\Vol{111}
\Year{1989}
\Pages{381--385}

\nextref
Oo74
\auth
A.I.,Ooms;
\endauth
\paper{On Lie algebras having a primitive universal enveloping algebra}
\journal{J. Algebra}
\Vol{32}
\Year{1974}
\Pages{488--500}

\nextref
Oo76
\auth
A.I.,Ooms;
\endauth
\paper{On Lie algebras with primitive envelopes, supplements}
\journal{Proc. Amer. Math. Soc.}
\Vol{58}
\Year{1976}
\Pages{67--72}

\nextref
Oo80
\auth
A.I.,Ooms;
\endauth
\paper{On Frobenius Lie algebras}
\journal{Comm. Algebra}
\Vol{8}
\Year{1980}
\Pages{13--52}

\nextref
Pre99
\auth
A.,Premet;S.,Skryabin;
\endauth
\paper{Representations of restricted Lie algebras and families of associative $\cal L$-algebras}
\journal{J.~Reine Angew. Math.}
\Vol{507}
\Year{1999}
\Pages{189--218}

\nextref
Rad76
\auth
D.E.,Radford;
\endauth
\paper{The order of the antipode of a finite dimensional Hopf algebra is finite}
\journal{Amer. J. Math.}
\Vol{98}
\Year{1976}
\Pages{333--355}

\nextref
Rad77
\auth
D.E.,Radford;
\endauth
\paper{Operators on Hopf algebras}
\journal{Amer. J. Math.}
\Vol{99}
\Year{1977}
\Pages{139--158}

\nextref
Rad92
\auth
D.E.,Radford;
\endauth
\paper{On the antipode of a quasitriangular Hopf algebra}
\journal{J. Algebra}
\Vol{151}
\Year{1992}
\Pages{1--11}

\nextref
Rad93
\auth
D.E.,Radford;
\endauth
\paper{Minimal quasitriangular Hopf algebras}
\journal{J. Algebra}
\Vol{157}
\Year{1993}
\Pages{285--315}

\nextref
Scha96
\auth
P.,Schauenburg;
\endauth
\paper{Hopf bigalois extensions}
\journal{Comm. Algebra}
\Vol{24}
\Year{1996}
\Pages{3797--3825}

\nextref
Scha97
\auth
P.,Schauenburg;
\endauth
\paper{A bialgebra that admits a Hopf-Galois extension is a Hopf algebra}
\journal{Proc. Amer. Math. Soc.}
\Vol{125}
\Year{1997}
\Pages{83--85}

\nextref
Schn90
\auth
H.-J.,Schneider;
\endauth
\paper{Principal homogeneous spaces for arbitrary Hopf algebras}
\journal{Isr. J.~Math.}
\Vol{72}
\Year{1990}
\Pages{167--195}

\nextref
St78
\auth
H.,Strade;
\endauth
\paper{Darstellungen aufl\"osbarer Lie $p$-Algebren}
\journal{Math. Ann.}
\Vol{232}
\Year{1978}
\Pages{15--32}

\nextref
St
\auth
H.,Strade;R.,Farnsteiner;
\endauth
\book{Modular Lie Algebras and their Representations}
\bookseries{Marcel Dekker Textbooks and Monographs}
\Vol{116}
\publisher{Marcel Dekker, New York}
\Year{1988}

\nextref
Sw
\auth
M.E.,Sweedler;
\endauth
\book{Hopf Algebras}
\publisher{Benjamin, New York}
\Year{1969}

\nextref
Taft75
\auth
E.J.,Taft;R.L.,Wilson;
\endauth
\paper{Hopf algebras with nonsemisimple antipode}
\journal{Proc. Amer. Math. Soc.}
\Vol{49}
\Year{1975}
\Pages{269--276}

\nextref
Taft80
\auth
E.J.,Taft;R.L.,Wilson;
\endauth
\paper{There exist Hopf algebras with antipodes of arbitrary even order}
\journal{J. Algebra}
\Vol{62}
\Year{1980}
\Pages{283--291}

\nextref
Ul82
\auth
K.-H.,Ulbrich;
\endauth
\paper{Galoiserweiterungen von nicht-kommutativen Ringen}
\journal{Comm. Algebra}
\Vol{10}
\Year{1982}
\Pages{655--672}

\nextref
Oy94
\auth
F.,Van Oystaeyen;Y.,Zhang;
\endauth
\paper{Galois-type correspondences for Hopf-Galois extensions}
\journal{$K$-Theory}
\Vol{8}
\Year{1994}
\Pages{257--269}

\nextref
Wei71
\auth
B.Yu.,Weisfeiler;V.G.,Kac;
\endauth
\paper{On irreducible representations of Lie $p$-algebras}
\journal{Functional Anal. Appl.}
\Vol{5}
\Year{1971}
\Pages{111--117}

\endreferences
\bye